\documentclass{article}

\usepackage[utf8]{inputenc}
\usepackage{geometry}
\usepackage{graphicx}%
\usepackage{multirow}%
\usepackage{amsmath,amssymb,amsfonts}%
\usepackage{amsthm}%
\usepackage{mathrsfs}%
\usepackage[title]{appendix}%
\usepackage{xcolor}%
\usepackage{textcomp}%
\usepackage{manyfoot}%
\usepackage{booktabs}%
\usepackage[boxruled,vlined]{algorithm2e}
\usepackage{algpseudocode}%
\usepackage{listings}%
\usepackage{hyperref}%
\usepackage{caption}%
\usepackage{subcaption}%
\usepackage{tikz}%
\usepackage{url}
\usepackage{xcolor}
\usepackage{hyperref}
\hypersetup{
    colorlinks=true,
    linkcolor=blue,
    filecolor=magenta,      
    urlcolor=cyan,
    citecolor=blue
}
\usepackage{bm}
 
\newcommand{\bbr}{\mathbb{R}}

\providecommand{\abs}[1]{\lvert#1\rvert}
\providecommand{\norm}[1]{\lVert#1\rVert}
\newcommand{\la}{\langle}
\newcommand{\ra}{\rangle}
\newcommand{\ba}{\bar{a}}

\newcommand{\bigslant}[2]{{\raisebox{.1em}{$#1$}\left/\raisebox{-.1em}{$#2$}\right.}}
\renewcommand{\circ}{\diamond}

\newtheorem{theorem}{Theorem}
\newtheorem{definition}{Definition}
\newtheorem{lemma}{Lemma}
\newtheorem{corollary}{Corollary}
\newtheorem{proposition}{Proposition}
\newtheorem{remark}{Remark}
\newtheorem{example}{Example}

\title{On Covering Euclidean Space with Q-arrangements of Cones}

\author{Khalil Ghorbal and Christelle Kozaily} 

\date{INRIA, France}

\begin{document}
\maketitle

\abstract{
This paper is concerned with a covering problem of Euclidean space by a particular arrangement of cones that are not necessarily full and are allowed to overlap. 
The problem provides an equivalent geometric reformulation of the solvability of the linear complementarity problem defining the class of Q-matrices. 
Assuming feasibility, we rely on standard tools from convex geometry to study maximal connected uncovered regions, we term \emph{holes}.  
We then use our approach to fully characterize the problem for dimension $3$, regardless of degeneracy. 
We further provide, for $n \leq 3$, an algebraic characterization for the class of Q-matrices. That is, we show that, $M$ is a Q-matrix if and only if its entries belong to an explicit semi-algebraic set (in dimension $9$) where all the involved polynomials are subdeterminants of $M$.  
We showcase the usefulness of such a characterization by generating $3$-by-$3$ Q-matrices with specific interesting properties on the involved cones. 
}

\paragraph{AMS/MSC.} 90C33, 68W30, 03C10.

\section*{Introduction}
Given a vector $q \in \bbr^n$ and an $n \times n$ matrix $M$ over the reals, the \emph{linear complementarity problem}, LCP$(q,M)$, asks whether there exists a pair $w,z \in \bbr^n$ satisfying $w-M z = q$, $w,z\geq 0$, and $w.z=0$, where $w,z \geq 0$ means that $w$ and $z$ belong to $\bbr^n_+$, the nonnegative orthant of $\bbr^n$, and $w.z$ is the scalar product of $w$ and $z$ (cf. \cite{cottle2009linear}). 
When LCP$(q,M)$ admits a solution, it is said to be \emph{solvable}. 
When a solution satisfying only $w,z \geq 0$ exists (i.e. when dropping the scalar product requirement), LCP$(q,M)$ is said to be \emph{feasible}. 
Related to the solvability and feasibility concepts, several classes of matrices were defined in the literature. Three classes are in particular relevant to this work. 
When LCP$(q,M)$ is feasible for all $q$, $M$ is called an \emph{S-matrix}. When LCP$(q,M)$ is solvable for all $q$, $M$ is called a \emph{Q-matrix}. If furthermore such a solution is unique for all $q$, $M$ is called a \emph{P-matrix}. 

The solvability of linear complementarity problems was tackled from different angles. \\ 
In \cite{COTTLE198173}, a focus on the linear application represented by $M$ proved to be useful. Degree theory~\cite[Chapter 6]{pang1979onq} was exploited for certain classes of structured matrices (e.g.~\cite{Garcia1983}), and more recently a new sufficient condition was provided in~\cite{radonsCueto}.   
More in line with this work, \cite{MURTY197265} pointed out an insightful geometric interpretation for LCP$(q,M)$.  
Instead of fixing $q \in \bbr^n$ and solving for $w$ and $z$, one could instead fix the pair $(w,z)$ and solve for the vectors $q$. 
The constraints $w,z \geq 0$ and $w.z=0$ imply that, for each $i$, either $w_i$ or $z_i$ has to vanish and the remaining component has to be nonnegative making $q$ an element of a cone spanned by some columns of $-M$ and $I$ (the identity matrix in dimension $n$). 
In turn, asking for a solution for each $q \in \bbr^n$ becomes equivalent to asking whether $\bbr^n$ is covered by the union of these cones. 

In the late fifties, \cite{Samelson1958APT} characterized a partition of Euclidean space with $2^n$ full cones using a separation condition. Their characterization provides a geometric reformulation of what later became known as the class of P-matrices (with the separation being captured as the positiveness of the principal minors of $M$).  
Several serious attempts have been subsequently made to characterize the class of Q-matrices concisely and, despite the rich literature devoted to the problem, it remains open even for low dimensions, see e.g.~\cite{fredricksen1986finite}. 
\cite{KELLY1979175} proved that, while the set of non-degenerate Q-matrices is open for $n=3$, it is not for $n=4$. 
Some subclasses of Q-matrices with adequate structures were easier to tackle. The problem of recognizing a P-matrix was shown to be co-NP-complete by \cite{Coxson1994ThePP}. A larger class than P-matrices relies on oriented matroids realized by $I$ and $-M$ was also studied in~\cite{watson1974variational,morris1986oriented}. \cite{morris1988counterexamples} provided a (counter) example for $n=4$ showing that the signs of subdeterminants of $M$ alone are not enough to characterize Q-matrices. 
As of today, the cost of checking generic Q-matrices remains prohibitive in practice, see~\cite{AganagicCottle78,NaimanStone1998}, and the best known practical approach presented in \cite{DeLoeraMorris1999} uses secondary and universal polytopes and is limited to dimensions less than $10$. 
It's worth noting that these theoretical facts about the hardness of the problem are in contrast with the complexity of solving LCP$(q,M)$ for a \emph{fixed} $q$ and $M$. Indeed, while the general problem (over the integers) has been shown to be NP-complete by~\cite{NPLCP89}, there are specific practical instances for which one can go beyond thousands of variables (cf. e.g.~\cite{luigiVincenzo08}). 

The rest of the paper is organized as follows. After a formal introduction of the problem (Section~\ref{sec:local}), we investigate the relatively simpler feasibility problem to arrive at useful cones we term \emph{minimal} (Section~\ref{sec:gamma}). Section~\ref{sec:necConds} focuses on solvability where we revisit the original \emph{separation condition} by~\cite{Samelson1958APT} in Section~\ref{sec:sepsep} and show a similar necessary (but not sufficient) condition for Q-matrices. 
This condition insinuates a lead to refine the standard concept of separation using intersections of some specific cones instead of hyperplanes separating vectors (Sections~\ref{sec:dyadcover} and~\ref{sec:self}). 
These different ingredients are leveraged in Section~\ref{sec:n=3} to give a complete characterization of Q-matrices in dimension $3$ regardless of degeneracy. 
Finally, section~\ref{sec:surrounding} provides an equivalent algebraic characterization in terms of signs of the subdeterminants of the involved matrix. 

\textbf{Contributions.} 
    To the best of our knowledge, the following are novel results related to Q-matrices. 
    Proposition~\ref{prop:gammacovering} gives a necessary condition useful to characterize S-matrices, i.e. the feasibility problem, by triangulating the space. 
    Proposition~\ref{prop:separation} provides a necessary separation condition for the Q-covering problem (cf. Definition~\ref{def:qmatrix}). 
    Proposition~\ref{prop:a1a1pcov} gives an interesting covering property of particular relevant cones we term \emph{minimal cones}. 
    The remaining results concern only dimension $3$. 
    Theorem~\ref{thm:Gcovered} characterizes uncovered regions assuming feasibility. 
    Corollary~\ref{coro:char} strengthens~\cite[Theorem 4.7]{Garcia1983} by dropping the strong non-degeneracy assumption. 
    Theorem~\ref{thm:symbolic} gives necessary and sufficient conditions for the Q-covering problem. 
    The algorithms in Section~\ref{sec:surrounding} provide an explicit algebraic characterization for Q-matrices as sign conditions of the subdeterminants of the involved matrix (Theorem~\ref{thm:sub3n}). 
    Such a characterization turned out to be very convenient to generate Q-matrices with interesting properties (cf. Examples~\ref{ex:examplen3} and~\ref{ex:exampleMurty}).  
 
\section{Q-covering}\label{sec:local}
Let $g_1,\dotsc,g_m \in \bbr^n$. 
The polyhedral cone or simply \emph{cone} (resp. linear subspace) spanned by the vectors $g_i$ will be denoted by $\la g_1,\dotsc,g_m\ra$ (resp. $(g_1,\dotsc,g_m)$).    
A cone is said to be \emph{non-degenerate} if its generators are linearly independent (as vectors in $\bbr^n$). It is \emph{degenerate} otherwise. 
A non-degenerate cone is said to be \emph{full} (or \emph{simplicial}) if its generators form a basis of $\bbr^n$.\footnote{We warn the reader that, some authors, e.g.~\cite{cottle2009linear},  refer to simplicial complementary cones as non-degenerate.} 
When $C=\la g_1,\dotsc,g_m\ra=(g_1,\dotsc,g_m)$, $C$ is said to be \emph{flat}. 
Flatness and degeneracy should not be confused. 
A degenerate cone is not necessarily flat (e.g. the half line $\la g_1, g_1\ra$, with $g_1 \neq 0$). 
Recall that a cone is said to be \emph{non-pointed} if it contains both a nonzero vector and its opposite. It is \emph{pointed} otherwise. 
Thus, while a flat cone is necessarily non-pointed, the converse is not necessarily true: for instance, a (closed) half-plane is both non-pointed and non-flat. 
The relative interior of a cone $C$ will be denoted by $C^\circ$.\footnote{In classical textbooks, e.g.~\cite{rockafellar1997convex}, $\la g_1,\dotsc,g_m\ra$ is denoted by $\operatorname{cone}\{g_1,\dotsc,g_m\}$ and the relative interior of $C$ is denoted by $\operatorname{ri} C$.} 
One proves that $\la g_1,\dotsc,g_m\ra^\circ$ is the set spanned by the positive linear combination of $g_1,\dotsc,g_m$. 

%
Let $\{a_1,a'_1\},\dotsc,\{a_n,a'_n\}$ denote $n$ pairs, or \emph{dyads}, of vectors in $\bbr^n$. 
Let the matrices $A$ and $A'$ denote respectively $\begin{pmatrix}
    a_1 & \dotsb & a_n
\end{pmatrix}$ and $\begin{pmatrix}
    a'_1 & \dotsb & a'_n
\end{pmatrix}$.  
Consider the mapping 
\begin{alignat*}{2}
[A,A']:\quad \{0,1\}^n &\to \bbr^{n \times n} \\
b &\mapsto [A,A']_b
\end{alignat*}
where the $i$th column vector of the matrix $[A,A']_b$ is $a_i$ if $b_i=1$ and $a'_i$ otherwise. 
A \emph{complementary cone}, or c-cone, $C$ is the 
cone spanned by the column vectors of $[A,A']_b$ for some valuation of $b$, that is $C = [A,A']_b(\bbr^n_+)$, the image set of $\bbr^n_+$ through the linear application represented by $[A,A']_b$. We call the column vectors of $[A,A']_b$, \emph{the generators} of the c-cone $C$. 
A \emph{complementary face}, or a c-face, is a face of a c-cone~\cite[Section 18]{rockafellar1997convex}. An element of a c-face is a nonnegative combination of $m$, $1 \leq m \leq n$, column vectors of $[A,A']_b$~\cite[Corollary 18.3.1]{rockafellar1997convex}. 
When the subspace spanned by a face $F$ has dimension $n-1$, $F$ is called a \emph{facet}. 
We define similarly a \emph{complementary linear subspace}, or c-subspace, as the set of linear combinations of $m$, $1 \leq m \leq n$, columns of $[A,A']_b$ for some $b \in \{0,1\}^n$. 

\begin{definition}[Covered and surrounded sets] 
A vector is said to be \emph{covered} if it belongs to a c-cone. 
A subset of $\bbr^n$ is said to be covered if all its vectors are covered. 
A vector is said to be \emph{surrounded} if it has a covered neighborhood. A subset of $\bbr^n$ is said to be surrounded if all its vectors are surrounded.
\end{definition}
In this work we investigate what conditions the pair $\{a_i,a'_i\}$ has to satisfy for Euclidean space to be covered. Formally, we are interested in the following problem.  
\begin{definition}[Q-covering]
\label{def:qmatrix}
Let $\{a_i,a'_i\}$, $i=1,\dotsc,n$ denote a list of $n$ dyads of vectors in $\bbr^n$. 
The \emph{Q-covering} decision problem asks whether $\bbr^n$ is covered, that is whether $\Sigma=\bbr^n$ where 
\[
\Sigma = \{a_1,a'_1\}\oplus \dotsb \oplus \{a_n,a'_n\} := \bigcup_{b \in \{0,1\}^n} [A,A']_b(\bbr^n_+) 
\enspace .
\]
\end{definition}
We adopt the sum notation for $\Sigma$ in the sequel. 
We observe that $\Sigma$ is invariant under any permutation of the indices ($\oplus$ is commutative). 
We think of the vector $a'_i$ as the \emph{dual} or the \emph{symmetric partner} of $a_i$ in the sense that $\Sigma$ remains invariant when swapping $a_i$ and $a'_i$ for any $i$.
The prime `$\ ^\prime$' can be thought of as an involutive operator providing the symmetric partner of $a_i$. 
For instance, the opposite is a very special `$\ ^\prime$' operator: when $a'_i = -a_i$, for all $i$, $\Sigma$ is the partition of $\bbr^n$ into the $2^n$ standard orthants. 

If $L$ is a non-singular matrix, then $q \in \bbr^n$ is covered if and only if $q$ belongs to a c-cone $\la a_1,\dotsc,a_n\ra$ say, which is equivalent to $L q \in \la L a_1,\dotsc, L a_n\ra$. Therefore, the Q-covering problem is invariant under non-singular linear transformations of the involved vectors. 
If all c-cones are degenerate, $\bbr^n$, as a Baire space, cannot be covered. 
So for $\Sigma$ to be covering, at least one c-cone must be full. 
There is thus no loss of generality in considering the standard basis $e_1,\dotsc,e_n$ of $\bbr^n$ respectively for $a_1,\dotsc,a_n$ (equivalently $A$ is the identity matrix) as required in the standard definition of the linear complementarity problem. Said differently, with respect to Definition~\ref{def:qmatrix}, $\{e_1,-M_1\}\oplus \dotsb \oplus \{e_n,-M_n\} = \bbr^n$ if and only if $M=(\begin{smallmatrix}
    M_1 & \dotsm & M_n
\end{smallmatrix})$ is a Q-matrix. 

\begin{remark}
Regardless of the exact geometric intersection between two c-cones, each c-cone has $n$ (abstract) \emph{neighbors} where two c-cones respectively generated by the column vectors of $[A,A']_{b}$ and $[A,A']_{b'}$ are neighbors if and only if the Hamming distance between $b$ and $b'$ is exactly one. 
Therefore, the c-cones can be put in correspondence with the vertices of an $n$-dimensional hypercube graph $Q_n$ where the neighboring relationship is represented by the adjacency of the vertices in $Q_n$. 
Stitching together all c-cones along their common abstract facets, one at a time, following their neighborhood relationship, amounts to following the longest Hamiltonian cycle of $Q_n$ (of length $2^n$). 
In general, a family of convex sets that covers the space and for which the aforementioned neighboring relation makes sense is said to form a \emph{Q-arrangement}. In~\cite{KELLY1979175}, Q-arrangements of full c-cones are studied.
\end{remark}

By introducing the equivalence relation $\simeq$ over $\bbr^n \setminus \{0\}$ defined by $u \simeq v$ if and only if $u = \lambda v$ for some positive scalar $\lambda$, one observes that if a nonzero vector $q$ is covered then so is $q' \simeq q$. 
Thus one can equivalently study the covering problem of the quotient space $\bigslant{(\bbr^n \setminus \{0\})}{\simeq}$ instead of $\bbr^n$ (observing that $q=0$ is trivially covered as it belongs to all c-cones). 
The quotient space $\bigslant{(\bbr^n \setminus \{0\})}{\simeq}$ is homeomorphic to the unit sphere $S_{n-1}$ where each half-line is represented by its unit generator. 
Similarly, each c-cone is represented by the (possibly degenerate) spherical $(n-1)$-simplex formed by the representatives of the $n$ generators of the cone. The spherical covering was for instance instrumental in~\cite{KELLY1979175,morris1988counterexamples,cottle1981spherically}. 
(Notice, however, that the collection of the so obtained complementary simplices does not necessarily form a simplicial complex, see~\cite{goerss2012simplicial}, since degeneracy and full dimensional intersections are allowed.) 

In the sequel, we will find it useful to fix one or several coordinates of the Boolean vector $b$. We define 
\[
\Sigma(a_i) :=  \bigcup_{\substack{b \in \{0,1\}^n \\ b_i=1}} [A,A']_b(\bbr^n_+),  \qquad \Sigma(a'_i) :=  \bigcup_{\substack{b \in \{0,1\}^n \\ b_i=0}} [A,A']_b(\bbr^n_+) \enspace .
\]
Clearly, for all $i$, $\Sigma = \Sigma(a_i) \cup \Sigma(a'_i)$. We say that $\Sigma(a_i)$ is the set of c-cones \emph{rooted} at $a_i$ to make the syntactic requirement of the definition salient. (Indeed, $a_i$ could be among the generators of a c-cone without necessarily having $b_i=1$. For instance when $a_2 \simeq a_1$, $a_1$ qualifies as a generator for $\la a'_1,a_2\ra$ while $b_1=0$.) 
The set of c-cones will be denoted by $\operatorname{cones}(\Sigma)$. Similarly, $\operatorname{cones}(\Sigma(a_i))$ will denote the set of c-cones rooted at $a_i$. 
We further let $\Sigma_k$, $1 \leq k \leq n$, denote the union of all c-faces with $m\leq k$ generators: for instance $\Sigma_1$ is the set of cones $\la a_i \ra$, $\la a'_i \ra$, $i=1,\dotsc,n$. 
By convention, we let $\Sigma_0$ denote the set of vectors $a_1,\dotsc,a_n$, $a'_1,\dotsc,a'_n$.  

\begin{definition}[Self surrounding]\label{def:self}
    The vector $a_i$ is said to be \emph{self surrounded} if it has a neighborhood covered by $\Sigma(a_i)$. 
\end{definition}
\begin{definition}[Lazy covering and surrounding]
    \label{def:lazy}
    We say that $a_i$ is \emph{lazily covered} if it belongs to $\Sigma(a'_i)$. 
    It is \emph{lazily surrounded} if it belongs to the topological interior of a c-cone rooted at $a'_i$. 
\end{definition}

We end this section by a simple definition which will be instrumental in the sequel. 
We use $\subseteq$ for set inclusion and $\subset$ for proper (or strict) set inclusion.   
\begin{definition}[Hole] 
A \emph{hole} is a non-empty open connected region in $\Sigma^c := \bbr^n \setminus \Sigma$, the complement of $\Sigma$ with respect to $\bbr^n$. 
A \emph{maximal} hole is a hole $K$ such that if $K'$ is another hole, $K \subseteq K'$ implies $K'=K$. 
\end{definition}
Since $\Sigma$ is a union of finitely many closed sets, $\Sigma^c$ is an open set. 
In general, $\Sigma^c$ is a union of disconnected maximal holes. In this work, we will consider only maximal holes and refer to them simply as holes. 
For instance, for $n=2$, $\Sigma = \{e_1,-e_2\} \oplus \{e_2,-e_1\}$, $\Sigma^c$ has two disconnected holes, namely $\la e_1, -e_2\ra^\circ$ and $\la -e_1, e_2\ra^\circ$. 
A hole is not necessarily convex. For instance, when $\Sigma = \{e_1,e_1\} \oplus \{e_2,e_2\} = \la e_1,e_2\ra$, $\Sigma^c = \bbr^2 \setminus \la e_1, e_2\ra$ is a non-convex hole.

\section{Feasibility}\label{sec:gamma} 
When $\bbr^n$ is covered by c-cones, some necessary conditions are intuitively clear and plausible. 
For instance, one can easily see that Euclidean space cannot be covered by $\Sigma$ when the $2^n$ c-cones are all not full or when all vectors of $\Sigma_0$ belong to the same half-space. 
These conditions among others, collectively insinuate a broader necessary condition requiring the vectors of $\Sigma_0$ to be ``well scattered'' in the space to form enough full c-cones that are in turn sufficiently distributed to achieve a covering. For instance, it is well known that at least $n+1$ full cones are required to cover $\bbr^n$: a full cone $C$ with an additional vector in the topological interior of $-C$ partition the space in $n+1$ full cones.~\footnote{This fact can be seen as a corollary of Stiemke's theorem \cite[Theorem 2.7.12]{cottle2009linear}.} We shall see that, in the context of this paper, a similar \emph{separation} property follows from Proposition~\ref{prop:separation}. 
Let $\Gamma = \la a_1,a'_1,\dotsc,a_n,a'_n\ra$. 
As we are working in dimension $n$, $\Gamma$ can be seen as the union of $\binom{2n}{n}$ cones generated by any $n$ vectors in $\Sigma_0$. 
The set of such cones will be denoted by $\operatorname{cones}(\Gamma)$. 
We proceed to investigate under which conditions $\Gamma$ covers $\bbr^n$ ensuring feasibility.
We do this in Proposition~\ref{prop:gammacovering} after stating a technical lemma (akin to fan and barycentric triangulation) characterizing flat cones as those having $0$ in their relative interior. 
\begin{lemma}\label{lem:flatcone}
    Suppose that $0 \in \la g_1,\dotsc,g_{m+1}\ra^\circ$, $1 \leq m \leq n$, then $\la g_1,\dotsc,g_{m+1}\ra$ is flat of dimension at most $m$. 
    Moreover, $\la g_1,\dotsc,g_{m+1}\ra = \cup_i G_i$ where $G_i$, $1 \leq i \leq m+1$, is the cone generated by $\{g_1,\dotsc,g_{m+1}\} \setminus \{g_i\}$. 
\end{lemma}
\begin{proof}
    There exist $\alpha_1,\dotsc,\alpha_{m+1}>0$ such that $0 = \sum_{i=1}^{m+1} \alpha_i g_i$. 
    Thus $(g_1,\dotsc,g_{m+1}) = (g_1,\dotsc,g_{m})$. 
    Let $x \in (g_1,\dotsc,g_{m})$. Then $x = \sum_{i=1}^m \lambda_i g_i$ where $\lambda_i \in \bbr$. 
    If $\lambda_i \geq 0$ for all $1 \leq i \leq m$, then $x \in \la g_1,\dotsc,g_m\ra \subseteq \la g_1,\dotsc,g_{m+1}\ra$. Otherwise, there exists a non-empty set of indices $J$ such that $\lambda_j < 0$ for all $j \in J$. 
    Let $\lambda = \max_{j \in J} \{\tfrac{-\lambda_j}{\alpha_j}\} > 0$. 
    So, for all $i$, $\lambda_i + \lambda\alpha_i \geq 0$. Then 
    \[
    x = x + \lambda\times 0 = \sum_{i=1}^m \lambda_i g_i + \lambda \sum_{i=1}^{m+1} \alpha_i g_i = \lambda\alpha_{m+1} g_{m+1} + \sum_{i=1}^m (\lambda_i + \lambda\alpha_i) g_i
    \]
    and $x \in \la g_1,\dotsc,g_{m+1}\ra$. Therefore $(g_1,\dotsc,g_{m+1}) = (g_1,\dotsc,g_{m}) \subseteq \la g_1,\dotsc,g_{m+1}\ra \subseteq (g_1,\dotsc,g_{m+1})$. This proves the first part of the lemma. 
    For the second part, if $J$ is empty, then $x \in G_{m+1}$, otherwise, if we let $k\in J$ denotes the index for which $\lambda = \tfrac{-\lambda_k}{\alpha_k}$, then one sees that $x \in G_k$. Thus $\la g_1,\dotsc,g_{m+1}\ra \subseteq G_{m+1} \cup (\cup_{j \in J} G_j) \subseteq \cup_{i=1}^{m+1} G_i$. The converse inclusion is immediate and the equality holds as stated.  
\end{proof}

\begin{proposition}\label{prop:gammacovering}
    Let $k \geq n+1$ and $g_1,\dotsc,g_k$ denote $k$ nonzero vectors of $\bbr^n$. 
    If $\bbr^n = \la g_1,\dotsc,g_k\ra$ then there exist $i_1, \dotsc, i_{m+1}$ with $1 \leq m \leq n$ such that $\la g_{i_1},\dotsc,g_{i_{m+1}} \ra$ is a flat of dimension $m$. 
\end{proposition}
\begin{proof}
    The proof is by induction on $n$. 
    For $n=1$, suppose $\bbr = \la g_1,\dotsc,g_k \ra$ with $k \geq 2$. 
    Then $0 \in \la g_1,\dotsc,g_k\ra^\circ$ and there exist $\lambda_1,\dotsc,\lambda_k > 0$ such that $\sum_i \lambda_i g_i = 0$. 
    If $g_i < 0$ for all $i$ we get a contradiction, so there exists $g_j > 0$. 
    Similarly, if $g_i > 0$ for all $i$ we also get a contradiction, so there exists $g_\ell < 0$ and we have $\bbr = \la g_j, g_\ell\ra$ and $m = n = 1$. 
    Let $\pi$ denote the orthogonal projection onto the hyperplane $g_1^\perp$ ($g_1 \neq 0$ by hypothesis). 
    Let $x \in g_1^\perp$ which is a subset of $\bbr^n = \la g_1,\dotsc,g_k\ra$. 
    Then $x = \sum_{i=1}^k \lambda_i g_i$ with $\lambda_i \geq 0$. 
    Thus $x = \pi(x) = \sum_{i=2}^k \lambda_i \pi(g_i)$ and therefore $g_1^\perp \subseteq \la \pi(g_2),\dotsc,\pi(g_k)\ra$. The converse inclusion is immediate by definition of $\pi$. Thus $g_1^\perp = \la \pi(g_2),\dotsc,\pi(g_k)\ra$. 
    By the induction hypothesis, there exists $1 \leq m \leq n-1$ such that $\la \pi(g_{i_2}),\dotsc,\pi(g_{i_{m+2}})\ra$ is a flat of dimension $m$. In particular $\pi(g_{i_2}),\dotsc,\pi(g_{i_{m+1}})$ are linearly independent. 
    Assume without loss of generality that $i_j=j$. 
    We have $0 = \sum_{j=2}^{m+2} \alpha_j \pi(g_j)$, with $\alpha_j > 0$. 
    Moreover, for each $g_i$, there exists $\theta_i \in \bbr$ such that $g_i = \theta_i g_1 + \pi(g_i)$. 
    Thus $0 = \sum_{j=2}^{m+2} \alpha_j (g_j - \theta_j g_1) = (-\sum_{j=2}^{m+2} \alpha_j\theta_j)g_1 + \sum_{j=2}^{m+2} \alpha_j g_j$. 
    Let $\gamma = -\sum_{j=2}^{m+2} \alpha_j\theta_j$. If $\gamma > 0$, then by Lemma~\ref{lem:flatcone}, $(g_1,\dotsc,g_{m+2}) = (g_1,\dotsc,g_{m+1}) =  \la g_1,\dotsc,g_{m+2}\ra$. 
    We further prove that $g_1,\dotsc,g_{m+1}$ are linearly independent. 
    Suppose $0 = \sum_{i=1}^{m+1} \sigma_i g_i$ then $0 = \sum_{i=2}^{m+1} \sigma_i \pi(g_i)$ and since $\pi(g_{2}),\dotsc,\pi(g_{{m+1}})$ are linearly independent, $\sigma_i=0$ for all $ 2 \leq i \leq m+1$. Thus $0 = \sigma_1 g_1$, $\sigma_1 =0$, and $\la g_1,\dotsc,g_{m+2}\ra$ is a flat of dimension $m+1$, which is what we wanted to prove. 
    If $\gamma = 0$, then by Lemma~\ref{lem:flatcone}, $(g_2,\dotsc,g_{m+1})=(g_2,\dotsc,g_{m+2}) = \la g_2,\dotsc,g_{m+2}\ra$. Thus $\la g_2,\dotsc,g_{m+2}\ra$ is a flat of dimension $m$. 
    If $\gamma < 0$, then $g_1 \in \la g_2,\dotsc,g_{m+2}\ra$ and $\bbr^n = \la g_2,\dotsc,g_k\ra$. 
    We repeat the same reasoning with $g_2$ and either we find a flat of dimension $m$ or $m+1$ or $g_2$ can be also removed from the list of generators till eventually reaching $n+1$ generators for $\bbr^n$ at which point $m=n$ and we are done. 
\end{proof}

Proposition~\ref{prop:gammacovering} will be used in Section~\ref{sec:surrounding} to characterize the feasibility in dimension $3$. 
In the sequel, assuming $\Gamma=\bbr^n$, we show in Proposition~\ref{prop:cara} that the space can be covered by what we term \emph{minimal cones}. 
Lemma~\ref{lem:usimeqv} is akin to the \emph{anti-exchange} property of $\bbr^n$~\cite[Chapter I, \S 3]{CoppelFoundationsOC} and Lemma \ref{lem:partition} can be seen as an application of the conical version of Caratheodory's theorem. 

\begin{lemma}\label{lem:usimeqv}
    Let $u,v,g_1,\dotsc,g_m \in \bbr^n$, $1 \leq m \leq n-1$, such that $u,g_1,\dotsc,g_m$ or $v,g_1,\dotsc,g_m$ are linearly independent. 
    Then $u \in \la v,g_1,\dotsc,g_m\ra$ and $v \in \la u,g_1,\dotsc,g_m\ra$ if and only if $u \simeq v$. When $m=1$, the linear independence condition can be dropped. 
\end{lemma}
\begin{proof}
    There exist $\alpha,\lambda, \alpha_i,\lambda_i \geq 0$ such that 
    $u = \alpha v + \sum_{i=1}^m \alpha_i g_i$ and $v = \lambda u + \sum_{i=1}^m \lambda_i g_i$. Assume $v,g_1,\dotsc,g_m$ are linearly independent (otherwise, swap $u$ and $v$ in what follows). 
    We get $(1 - \lambda\alpha)v = \sum_{i=1}^m(\lambda \alpha_i + \lambda_i) g_i$. 
    The linear independence implies $1- \alpha\lambda = 0$ and $\alpha \lambda_i + \alpha_i = 0$ for all $i$. Thus $\alpha,\lambda > 0$ and $\alpha_i=\lambda_i=0$ for all $i$ proving the statement. 
    For $m=1$, dropping the linear independence hypothesis, the proof is by case distinction. If $\lambda\alpha_1 + \lambda_1 = 0$ then $\alpha,\lambda > 0$ and $\alpha_1 = \lambda_1 = 0$ and therefore $v \simeq u$. Likewise, if $\alpha\lambda_1 +\alpha_1 = 0$, $u \simeq v$. 
    Suppose $\lambda\alpha_1 + \lambda_1 > 0$ and $\alpha\lambda_1 +\alpha_1 > 0$. If  $1 - \lambda \alpha < 0$ then $u \simeq -g_1$ and $v \simeq -g_1$ so $u \simeq v$. Otherwise $u \simeq v \simeq g_1$.
\end{proof}

\begin{lemma}\label{lem:partition}
    Let $G = \la g_1,\dotsc,g_n\ra$ denote a cone in $\bbr^n$. 
    Let $g$ be a nonzero vector in $G$ and let $G_i$ denote the cone generated by $\{g,g_1,\dotsc,g_n\} \setminus \{g_i\}$, $1 \leq i \leq n$. Then $G = \cup_{i} G_i$. 
    If moreover $G$ is full then there exists a non-empty set of indices $J \subseteq \{1,\dotsc,n\}$ such that $G = \cup_{j \in J} G_j$, $G_j$ is full, and $G_j \subseteq G$, for all $j \in J$. Moreover $G_j = G$ if and only if $g \simeq g_j$. 
\end{lemma}
\begin{proof}
    We have $g=\sum_{i=1}^n \alpha_i g_i$ for some nonnegative coefficients $\alpha_i$ and, since $g \neq 0$, there exists a non-empty set of indices $J \subseteq \{1,\dotsc,n\}$ such that $\alpha_j > 0$ for all $j \in J$. 
    Let $x \in G$, i.e. $x=\sum_{i=1}^n \beta_i g_i$ for some nonnegative coefficients $\beta_i$. 
    Let $\lambda = \min_{j \in J} \bigl\{\tfrac{\beta_j}{\alpha_j}\bigr\} = \tfrac{\beta_{k}}{\alpha_{k}}$ for some $k \in J$. We then have, for all $i=1,\dotsc,n$, $\beta_i -\lambda \alpha_i \geq 0$ and 
    \[
    x = \sum_{i=1}^n \beta_i g_i = \lambda g + \sum_{i=1}^n (\beta_i -\lambda \alpha_i) g_i = 
    \lambda g + \sum_{i\neq k} (\beta_i -\lambda \alpha_i) g_i \enspace .
    \]
    Thus $x \in G_k$ and $G \subseteq \cup_{j \in J} G_j \subseteq \cup_i G_i$. This proves the first statement. Suppose $G$ is full and $G_j$ is degenerate for an index $j \in J$, then $g$ belongs to the hyperplane $H_j$ generated by $g_i$, $i \neq j$.~\footnote{In this paper, a \emph{hyperplane} denotes a linear (or vector) subspace of dimension $n-1$.} But then $g_j \in H_j$ (because $\alpha_j > 0$) and $G$ is itself degenerate, a contradiction. 
    Thus $G_j$ is full for all $j \in J$. 
    Since $g \in G$, $G_j \subseteq G$ for all $j$. 
    If in addition $G \subseteq G_j$ for some index $j$, then by Lemma~\ref{lem:usimeqv} (applied to $g$ and $g_j$ as $u$ and $v$), this is equivalent to $g \simeq g_j$. 
\end{proof}

\begin{proposition}\label{prop:cara}
    Assume $\Gamma = \bbr^n$. Then for any $x \in \bbr^n$, there exists a \emph{minimal cone} $G \in \operatorname{cones}(\Gamma)$ containing $x$, that is $G$ is full and for any other full cone $G' \in \operatorname{cones}(\Gamma)$, $G' \subseteq G$ implies $G'=G$. 
\end{proposition}
\begin{proof}
    Since $\Gamma = \bbr^n$, by the Caratheodory theorem, there exists a full cone $G_1$ with $n$ generators in $\Sigma_0$ such that $x \in G_1$ ($G_1$ needs not be unique). 
    Let $g_1,\dotsc,g_n \in \Sigma_0$ denote the $n$ generators of $G_1$ and suppose there exists a nonzero vector $g_{n+1} \in \Sigma_0 \setminus \{g_1,\dotsc,g_n\}$ such that $g_{n+1} \in G_1$ and $g_{n+1} \not\simeq g_i$ for all $1 \leq i \leq n$. 
    By Lemma~\ref{lem:partition}, there exists a full cone $G_2 \subset G_1$ containing $x$. 
    Suppose without loss of generality that $g_1 \not\in G_2$ (such a vector must exist by construction of $G_2$), so $G_2 = \la g_2,\dotsc,g_{n+1}\ra$. 
    If $G_2$ itself contains a vector $g_{n+2} \in \Sigma_0 \setminus \{g_1,\dotsc,g_{n+1}\}$ such that $g_{n+2} \not\simeq g_i$ for all $2 \leq i\leq n+1$, then again by Lemma~\ref{lem:partition}, there exists a full cone $G_3 \subset G_2$ containing $x$. 
    At step $k$, we ask if $G_k$ contains a vector $g_{n+k} \in \Sigma_0 \setminus \{g_1,\dotsc,g_{n+k-1}\}$ such that $g_{n+k} \not\simeq g_i$ for all $k \leq i \leq n+k-1$.  
    The process must terminate after at most $n$ steps since $\Sigma_0 \setminus \{g_1,\dotsc,g_{n+n}\}$ is empty. 
    We let $G=G_k$, $1 \leq k \leq n$, where $G_k$ does not contain any vector from $\Sigma_0$ (apart from those in the equivalence classes of its generators with respect to $\simeq$).  
    In particular there is no other cone $G' \in \operatorname{cones}(\Gamma)$ such that $G' \subset G$. 
\end{proof}

Minimal cones play an important role in locating holes. They will be for instance instrumental in Section~\ref{sec:n=3} to characterize Q-covering for $n=3$. 

\begin{proposition}\label{prop:cara1}
    Assume $\Gamma = \bbr^n$. 
    If $\Sigma \subset \bbr^n$ then there exists a minimal cone $G \in \operatorname{cones}(\Gamma) \setminus \operatorname{cones}(\Sigma)$ such that $G \cap \Sigma^c$ is non-empty. 
\end{proposition}
\begin{proof}
    Let $x \in \Sigma^c$. By Proposition~\ref{prop:cara}, there exists a minimal cone $G \in \operatorname{cones}(\Gamma)$ containing $x$. The cone $G$ cannot be a c-cone since it contains $x$. So $G \not\in \operatorname{cones}(\Sigma)$. Clearly we have $x \in G \cap \Sigma^c$ and the intersection is non-empty. 
\end{proof}

\begin{remark}
    One can tighten the statement of Proposition~\ref{prop:cara1} by saying that $G^\circ \cap \Sigma^c$ is a non-empty open set. 
    If $x$ belongs to the interior of $G$ then $G^\circ \cap \Sigma^c$ is non-empty. Otherwise, since $\Sigma^c$ is an open set, a small perturbation of $x$ remains in $\Sigma^c$ while avoiding the boundaries of $G$. Thus $G^\circ \cap \Sigma^c$ is also non-empty. 
\end{remark}


While $\Gamma$ and $\Sigma$ are both finite unions of closed convex cones, they are only seemingly similar. An important difference between the two being convexity: $\Gamma$ is convex by definition, but $\Sigma$ is not. 
This difference introduces a substantial complexity to the covering problem. For instance, the extreme rays of $\Gamma$ (if any) are necessarily among the generators of $\Gamma$~\cite[Corollary 18.3.1]{rockafellar1997convex}. This same property is far from obvious for $\Sigma$ as we shall see in Section~\ref{sec:n=3}. 
The next section investigates some interesting properties of $\Sigma$.

\section{Solvability}\label{sec:necConds} 

We start by proving a necessary condition for $\bbr^n$ to be covered, namely, that for all $i$, $a_i,a'_i \neq 0$ and $a'_i \not\simeq a_i$ (geometrically, this means that $\la a_i,a'_i \ra$ is not reduced to the origin nor is it a half-line). 
\begin{lemma}
\label{lem:inseparable}
    Assume there exists at least one full c-cone. If all full c-cones meet at a nonzero vector then $\Sigma \subset \bbr^n$. 
    Moreover, if there exists an index $i$ such that $a_i=0$ or $a'_i \simeq a_i$ then $\Sigma \subset \bbr^n$. 
\end{lemma}
\begin{proof}
    If $a_i=a'_i=0$ for some index $i$, then all c-cones are degenerate contradicting the assumption. 
    Assume next that, there exists at least one full c-cone and that for all $i$, either $a_i\neq 0$ or $a'_i \neq 0$. 
    Suppose there exists a vector $p \neq 0$ that belongs to all full c-cones. The full c-cones cannot cover $-p$ as this would contradict the non-degeneracy assumption. This means that an open neighborhood $U$ of $-p$ is left to be covered by a finite number of degenerate cones, which is  impossible ($\bbr^n$ is a Baire space).~\footnote{Every complete metric space (such as Euclidean space) is a Baire space in which countable unions of closed sets with empty interior have also an empty interior.}
    In particular, if $a_i=0$ and $a'_i \neq 0$, then all full c-cones meet at $p \simeq a'_i$ and the covering is impossible. 
    Moreover, when $a_i \neq 0$ and $a'_i \simeq a_i$ then all c-cones meet at $p \simeq a_i$ and the covering is also impossible. 
\end{proof} 

Proposition~\ref{prop:separation} below states a necessary condition on pairs $\{a_i,a'_i\}$ 
drawing upon an existing result by~\cite{Samelson1958APT} characterizing the partition of Euclidean space into c-cones using \emph{separation} which we now define. 

\begin{definition}[Separation]
    Let $H$ denote a hyperplane and let $h$ denote a nonzero vector orthogonal to $H$. Two vectors $u,v$ are said to be \emph{separated} by $H$ if and only if the scalar products $u.h$ and $v.h$ are not zero and have opposite signs. (Geometrically, $u$ and $v$ belong each to a distinct open half-space bounded by $H$.) 
\end{definition}

\subsection{A Necessary Separation Condition}\label{sec:sepsep}

\begin{proposition}\label{prop:separation}
    If $\Sigma = \bbr^n$ then for each index $i=1,\dotsc,n$ there must exist a complementary hyperplane (that is a c-subspace of dimension $n-1$) separating the pair $\{a_i,a'_i\}$.  
\end{proposition}
\begin{proof} 
    By Lemma~\ref{lem:inseparable}, for all $i$, $a_i,a'_i \neq 0$ and $a_i \not\simeq a'_i$.  
    The proof is by contradiction. 
    Fix $i$ and let $V_r$, $r=1,\dotsc,2^{n-1}$, denote the c-subspaces generated by the remaining $(n-1)$ vectors $a_j,a'_j$, $j \neq i$. 
    If $\dim(V_r) < n-1$ for all $r$, then all c-cones are degenerate and $\Sigma \subset \bbr^n$. 
    Thus, there must exist a non-empty subset of indices $S \subseteq \{1,\dotsc,2^{n-1}\}$ such that $V_s$ is a hyperplane for all $s \in S$. 
    Assume that the pair $\{a_i,a'_i\}$ is not separated by any hyperplane $V_s$, $s \in S$. 
    If $a_i,a'_i \in V_s$ for all $s$, then again, all c-cones are degenerate and $\Sigma \subset \bbr^n$. 
    Thus, there must exist a set of indices $K \subseteq S$ such that  $a_i \not\in V_k$ or $a'_i \not\in V_k$ for all $k \in K$. 
    Let $V_k(i)$, $k \in K$, denote the closed half-space bounded by $V_k$ and containing either $a_i$ or $a'_i$ in its interior. 
    We thus have $a_i,a'_i \in V:=\cap_{k \in K} V_k(i)$. 
    The boundary of $V$, denoted $\partial V$, is a subset of the union of the boundaries of $V_k(i)$. 
    Since $V$ is a (closed) convex cone then $\la a_i,a'_i\ra \subseteq V$.   
    We show by contradiction that the topological interior of $V$ is non-empty. 
    Assume it is empty. Then $V$ is equal to its boundary $\partial V$ and $V = \partial V \subseteq \cup_{k \in K} V_k$. 
    Since $a'_i \in V$, there exists an index $k_1 \in K$ such that $a'_i \in V_{k_1}$ and $a_i \not\in V_{k_1}$. Thus $\dim(a_i,a'_i)=2$  
    and we construct a sequence of vectors $p_j \neq 0$, $1 \leq j \leq \abs{K}$, such that $p_j \in \la a_i,p_{j-1}\ra^\circ$ with $p_1 = a'_i$. 
    For each $j$, $p_j \in V$, so there exists an index $k_j \in K$ such that $p_j \in V_{k_j}$. 
    If $k_j = k_\ell$ for $1 \leq j,\ell \leq \abs{K}$, $j \neq \ell$, then $p_j,p_\ell \in V_{k_j} = V_{k_\ell}$, and $(a_i,a'_i)=(p_j,p_{\ell}) \subseteq V_{k_j}$. But then $a_i,a'_i \in V_{k_j}$, contradicting the fact that $k_j \in K$. 
    Thus $k_1,\dotsc,k_{\abs{K}}$ are all distinct making $K = \{k_1,\dotsc,k_{\abs{K}}\}$. 
    Now consider a vector $p \neq 0$ in $\la a_i,p_{\abs{K}}\ra^\circ$. 
    Then there exists $k_j \in K$ such that $p \in V_{k_j}$ and therefore $(a_i,a'_i)=(p,p_j) \subseteq V_{k_j}$, a contradiction. 
    So the topological interior of $V$ is non-empty. We denote it by $V^\circ$.  
    We show next that $-V^\circ$ is a hole. 
    First, by definition of $V^\circ$, $V^\circ \subseteq V_k(i)$, thus $-V^\circ \subseteq V_k(i)^c$, for each $k \in K$. Therefore $-V^\circ \subseteq \cap_k V_k(i)^c =  (\cup_k V_k(i))^c$. 
    Second, let $C(a_i,V_r)$ denote the c-cone spanned by $a_i$ and the generators of $V_r$. For $k \in K$, since $V_k$ does not separate the pair $\{a_i,a'_i\}$ then $C(a_i,V_k) \cup C(a'_i,V_k) \subseteq V_k(i)$. Thus 
    \[
    \Sigma_K := \cup_{k \in K} \bigl(C(a_i,V_k) \cup C(a'_i,V_k)\bigr) \subseteq \cup_{k \in K} V_k(i),  
    \]
    and $\bigl( \cup_{k \in K} V_k(i) \bigr)^c \subseteq (\Sigma_K)^c$. Thus $-V^\circ \subseteq (\Sigma_K)^c$. 
    Finally, all is left to cover $-V^\circ$ is a finite set of degenerate c-cones, namely $C(a_i,V_r)$ and $C(a'_i,V_r)$, $r\not\in K$. Covering an open set with a finite number of degenerate cones is impossible (by Baire's theorem). So there must exist an index $s \in S$ such that $V_s$ is a hyperplane separating the pair $\{a_i,a'_i\}$. 
\end{proof}

With respect to the notations of the proof of Proposition~\ref{prop:separation}, the seminal work of~\cite{Samelson1958APT} shows that $\Sigma$ partitions $\bbr^n$ if and only if, for each index $i=1,\dotsc,n$ and for each $r=1,\dotsc,2^{n-1}$, the c-subspace $V_r$ is a hyperplane separating the pair $\{a_i,a'_i\}$. 
Proposition~\ref{prop:separation} shows that, for $\Sigma$ to be (only) covering, it is necessary that for each index $i$ at least one c-subspace $V_r$ is a complementary hyperplane separating the pair $\{a_i,a'_i\}$. 
We know that this necessary condition is not sufficient in general to prove that $\Sigma$ is covering. For instance, in the plane, let $a'_1 = -e_1+e_2$ and $a'_2=-e_1$. 
Then $\Sigma = \{e_1,a'_1\}\oplus\{e_2,a'_2\}$ covers only the closed upper half-plane. 
The pair $\{e_1,a'_1\}$ is separated by the line $(e_2)$ and the pair $\{e_2,a'_2\}$ is separated by the line $(a'_1)$. 

In the next section we shift focus on how c-facets (instead of c-hyperplanes) might ``separate'' pairs via intersections. 

\subsection{Dyadic Covering}\label{sec:dyadcover}
We shall see several interesting properties concerned with intersections of c-facets with the cone $\la a_i, a'_i\ra$. 
We start by stating three generic relevant intersections (all depicted in Fig.~\ref{fig:intersections}) that will be instrumental in the rest of the paper. Lemmas~\ref{lem:front} and \ref{lem:side} play a role in the forthcoming Propositions~\ref{prop:a1a1pcov} and \ref{lem:a1ap1}.  Lemma~\ref{lem:back} shall be used in Section~\ref{sec:n=3}. 

\begin{figure}[bt]
     \centering
     \begin{subfigure}[b]{0.3\textwidth}
         \centering
         \begin{tikzpicture}
\draw[thick,-stealth] (0,0,0)--(1,0,0) node[right]{$v'_1$};  
\draw[-stealth] (0,0,0)--(-1,0,0);
\draw[thick,-stealth] (0,0,0)--(0,0,1) node[left]{$v_1$}; 
\draw[-stealth] (0,0,0)--(0,0,-1); 
\draw[thick,-stealth] (0,0,0)--(0,1,0) node[left]{$v_2$}; 
\draw[thick,-stealth] (0,0,0)--(1/2,-1/4,1/2) node[right]{$v_3$}; 
\coordinate[] (O) at (0,0,0);
\coordinate[] (api) at (1,0,0);
\coordinate[] (ai) at (0,0,1);
\coordinate[] (aj) at (0,1,0);
\coordinate[] (ak) at (1/2,-1/4,1/2);
\draw[fill=orange,opacity=0.3] (O)--(aj)--(ak);
\draw[fill=blue,opacity=0.3] (O)--(ai)--(api);
\end{tikzpicture}
         \caption{Front}
         \label{fig:front}
     \end{subfigure}
     \hfill
     \begin{subfigure}[b]{0.3\textwidth}
         \centering
         \begin{tikzpicture}
\draw[thick,-stealth] (0,0,0)--(1,0,0) node[right]{$v'_1$};  
\draw[-stealth] (0,0,0)--(-1,0,0);
\draw[-stealth] (0,0,0)--(0,0,1); 
\draw[thick,-stealth] (0,0,0)--(0,0,-1) node[right]{$-v_1$}; 
\draw[thick,-stealth] (0,0,0)--(0,1,0) node[left]{$v_2$}; 
\draw[thick,-stealth] (0,0,0)--(1/2,-1/2,-1/2) node[left]{$v_3$}; 
\coordinate[] (O) at (0,0,0);
\coordinate[] (api) at (1,0,0);
\coordinate[] (mai) at (0,0,-1);
\coordinate[] (aj) at (0,1,0);
\coordinate[] (ak) at (1/2,-1/2,-1/2);
\draw[fill=orange,opacity=0.3] (O)--(aj)--(ak);
\draw[fill=blue,opacity=0.3] (O)--(mai)--(api);
\end{tikzpicture}
         \caption{Side}
         \label{fig:side}
     \end{subfigure}
     \hfill
     \begin{subfigure}[b]{0.3\textwidth}
         \centering
         \begin{tikzpicture}
\draw[-stealth] (0,0,0)--(1,0,0);  
\draw[thick,-stealth] (0,0,0)--(-1,0,0) node[left]{$-v'_1$};
\draw[-stealth] (0,0,0)--(0,0,1); 
\draw[thick,-stealth] (0,0,0)--(0,0,-1) node[right]{$-v_1$}; 
\draw[thick,-stealth] (0,0,0)--(0,1,0) node[left]{$v_2$}; 
\draw[thick,-stealth] (0,0,0)--(-1,-.7,-1) node[left]{$v_3$}; 
\coordinate[] (O) at (0,0,0);
\coordinate[] (api) at (1,0,0);
\coordinate[] (mapi) at (-1,0,0);
\coordinate[] (mai) at (0,0,-1);
\coordinate[] (aj) at (0,1,0);
\coordinate[] (ak) at (-1,-.7,-1);
\draw[fill=orange,opacity=0.3] (O)--(aj)--(ak);
\draw[fill=blue,opacity=0.3] (O)--(mai)--(mapi);
\end{tikzpicture}
         \caption{Back}
         \label{fig:back}
     \end{subfigure}
        \caption{Conical intersections}
        \label{fig:intersections}
\end{figure}


\begin{lemma}[Front Intersection (cf. Fig~\ref{fig:front})]\label{lem:front}
    Fix $m \geq 2$ and let $C = \la v_1,v_2,\dotsc,v_m\ra$ and $C' = \la v'_1,v_2,\dotsc,v_m\ra$. 
    If there exists a nonzero vector in $\la v'_1,v_1 \ra \cap \la v_2,\dotsc,v_m\ra$ then $\la v'_1,v_1,v_2,\dotsc,v_{m}\ra = C \cup C'$. 
\end{lemma}
\begin{proof}
    Let $c \in \la v'_1,v_1,v_2,\dotsc,v_m \ra$. 
    There exist $\alpha'_1,\alpha_1,\dotsc,\alpha_m \geq 0$ such that $c = \alpha_1 v_1 + \alpha'_1 v'_1 + \sum_{i \geq 2} \alpha_i v_i$.  
    By assumption, there exist $\lambda'_1,\lambda_1,\dotsc,\lambda_m \geq 0$, $(\lambda_1,\lambda'_1) \neq (0,0)$, such that $\lambda_1 v_1 +\lambda'_1 v'_1 = \sum_{i \geq 2} \lambda_i v_i$. 
    Suppose $\lambda_1,\lambda'_1 > 0$ and let $\beta = \min\{\tfrac{\alpha_1}{\lambda_1},\tfrac{\alpha'_1}{\lambda'_1}\} \geq 0$. Then 
    \begin{align*}  
    c = c + \beta \times 0 &= \alpha_1 v_1 +\alpha'_1 v'_1 + \sum_{i \geq 2} \alpha_i v_i + \beta (-\lambda_1 v_1 -\lambda'_1 v'_1 + \sum_{i \geq 2} \lambda_i v_i) \\
    &=(\alpha_1 - \beta \lambda_1) v_1 + (\alpha'_1 - \beta \lambda'_1) v'_1 + \sum_{i \geq 2} (\alpha_i + \beta \lambda_i) v_i \enspace .
    \end{align*}
    By definition of $\beta$, the coefficients of $v_1$ and $v'_1$ are both nonnegative and one of them has to vanish. Thus $c \in C \cup C'$. 
    The same holds when $\lambda_1=0$ (resp. $\lambda'_1=0$) by taking $\beta$ as $\tfrac{\alpha'_1}{\lambda'_1}$ (resp. $\tfrac{\alpha_1}{\lambda_1}$). 
    The converse inclusion is immediate. 
\end{proof}


\begin{lemma}[Side Intersection (cf. Fig~\ref{fig:side})]\label{lem:side}
    Fix $m \geq 2$ and let $C = \la v_1,v_2,\dotsc,v_{m}\ra$ and $C' = \la v'_1,v_2,\dotsc,v_{m}\ra$. 
    If $\la -v_1,v'_1 \ra^\circ \cap \la v_2,\dotsc,v_{m}\ra$ is non-empty then $\la v_1,v'_1,v_2,\dotsc,v_{m}\ra \subseteq C$. 
\end{lemma}
\begin{proof}
    By assumption, there exist $\alpha'_1,\alpha_1,\dotsc,\alpha_{m} \geq 0$ such that $\alpha_1 (-v_1) +\alpha'_1 v'_1 =  \sum_{i \geq 2} \alpha_i v_i$ with  $\alpha_1,\alpha'_1 > 0$. 
    Thus $v'_1 \in C$ and $\la v_1,v'_1,v_2,\dotsc,v_{m} \ra \subseteq C$. 
\end{proof}

\begin{lemma}[Back Intersection (cf. Fig~\ref{fig:back})]\label{lem:back}
    Let $G = \la v_1,v'_1,v_2,\dotsc,v_{n-1}\ra$ denote a full cone in $\bbr^n$. 
    Suppose there exists a vector $p \neq 0$ such that $p \in \la -v_1,-v'_1 \ra \cap \la v_2,\dotsc,v_{n} \ra$ then $G^\circ \subseteq (\la v_1,v_2,\dotsc,v_n\ra \cup \la v'_1,v_2,\dotsc,v_n\ra)^c$. 
\end{lemma}
\begin{proof}
    There exist $\lambda_1,\lambda'_1,\lambda_2,\dotsc,\lambda_n \geq 0$, $(\lambda_1,\lambda'_1) \neq (0,0)$ such that 
    \begin{equation}\label{eq:bout1}
        p = \lambda_1 (-v_1) + \lambda'_1 (-v'_1) = \sum_{i=2}^n \lambda_i v_i \enspace .
    \end{equation}
    If $\lambda_n = 0$ then $G$ is degenerate. Thus $\lambda_n > 0$. 
    Let $g \in G^\circ$. Then $g = \alpha_1 v_1 + \alpha'_1 v'_1 + \sum_{i=2}^{n-1}\alpha_i v_i$ with $\alpha_1,\alpha'_1,\alpha_2,\dotsc,\alpha_{n-1} > 0$. 
    Suppose that $g \in \la v_1,\dotsc,v_n\ra$, then there exist $\beta_i \geq 0$ such that $g=\sum_{i=1}^n\beta_i v_i$. Thus $\alpha_1 v_1 + \alpha'_1 v'_1 + \sum_{i=2}^{n-1}\alpha_i v_i = \sum_{i=1}^n\beta_i v_i$ and 
    \begin{equation}\label{eq:bout2}
    \alpha'_1 v'_1 + \sum_{i=1}^{n-1}(\alpha_i - \beta_i) v_i = \beta_n v_n \enspace .
    \end{equation}
    If $\beta_n = 0$ then $\alpha'_1=0$ (because $v_1,v'_1,v_2,\dotsc,v_{n-1}$ are linearly independent), which contradicts $\alpha'_1 > 0$. Thus $\beta_n > 0$. 
    Since $\lambda_n,\beta_n > 0$ we can eliminate $v_n$ in Eq.~\ref{eq:bout1} and Eq.~\ref{eq:bout2}. We get 
    \[
    (\beta_n\lambda'_1 + \lambda_n\alpha'_1)(-v'_1) + \sum_{i=1}^{n-1}(\beta_n\lambda_i + \lambda_n(\alpha_i - \beta_i)) (-v_i) = 0 \enspace . 
    \]
    As $v_1,v'_1,v_2,\dotsc,v_{n-1}$ are linearly independent, all the coefficients must vanish. In particular,  $\beta_n\lambda'_1 + \lambda_n\alpha'_1 = 0$ which implies $\lambda_n\alpha'_1 = 0$, a contradiction. Thus $g \not\in \la v_1,\dotsc,v_n\ra$. 
    Swapping $v_1$ and $v'_1$ in the discussion above, we prove that $g \not\in \la v'_1,v_2,\dotsc,v_n\ra$ and therefore $G^\circ \subseteq (\la v_1,\dotsc,v_n\ra \cup \la v'_1,v_2,\dotsc,v_n\ra)^c$ as stated. 
\end{proof}

Proposition~\ref{prop:a1a1pcov} states an interesting covering property of the cone $\la a_i, a'_i\ra$. 
Its proof requires the following technical lemma.~\footnote{The original proof was by induction on $m$. The provided shorter constructive proof was suggested by Jean-Charles Gilbert in a private communication with the first author.}  
\begin{lemma}\label{lem:pathout3}
    Let $G = \la g_1,\dotsc,g_m\ra$, $m\geq 2$, denote a cone in $\bbr^n$ and let $q \not\in G$ be such that $\la g_1,q \ra$ is non-flat. 
    Suppose there exists a nonzero vector in $\la g_1,q \ra^\circ \cap G$ then there exists a nonzero vector in $\la g_2,\dotsc,g_m \ra \cap \la g_1,q \ra^\circ$.  
\end{lemma}
\begin{proof}
    Let $p' \neq 0$ be in $ \la g_1,q\ra^\circ \cap G $. Then there exists $\lambda_1,\lambda > 0$ and $\alpha_1,\dotsc,\alpha_m \geq 0$ such that $p' = \lambda_1 g_1 + \lambda q = \sum_{i=1}^m \alpha_i g_i$. 
    Since $q \not\in G$ then $0 < \lambda_1 -  \alpha_1$. 
    Let $p = (\lambda_1-\alpha_1) g_1 + \lambda q$. Then, $p \in \la g_2,\dotsc,g_m \ra \cap \la g_1,q \ra^\circ$. If $p=0$ then $-g_1 \simeq q$ ($q \neq 0$ because $q \not\in G$), making $\la g_1,q \ra$ flat. Thus $p \neq 0$. 
\end{proof}

\begin{proposition}\label{prop:a1a1pcov}
    The cone $\la a_i,a'_i\ra$ cannot be partially covered, i.e. either $\la a_i,a'_i\ra \subseteq \Sigma$ or $\la a_i,a'_i \ra^\circ \subseteq \Sigma^c$. 
\end{proposition}
\begin{proof}
    We fix $i$ to $1$ for clarity. 
    If $\dim(a_1, a'_1) \leq 1$, then $\la a_1, a'_1\ra \subseteq \Sigma$. 
    Assume next that $\dim(a_1,a'_1) = 2$ (in particular $a'_1 \not\simeq -a_1$). 
    Suppose there exists $a \in \la a_1, a'_1\ra^\circ$ such that $a \in \Sigma$. 
    So $a$ belongs to a c-cone $C=\la a_1,\dotsc,a_n\ra$ say (otherwise swap $a_1$ and $a'_1$ in what follows). 
    If $a'_1 \in C$ then, by convexity of $C$, $\la a_1, a'_1\ra \subseteq C \subseteq \Sigma$. 
    Assume $a'_1 \not\in C$.  
    By Lemma~\ref{lem:pathout3}, $\la a_1,a'_1\ra^\circ$ intersects the face $\la a_2,\dotsc,a_n\ra$ of $C$ in a nonzero vector. 
    By Lemma~\ref{lem:front}, $\la a'_1,a_1,\dotsc,a_n\ra \subseteq \la a_1,a_2,\dotsc,a_n \ra \cup \la a'_1,a_2,\dotsc,a_n\ra$ which implies $\la a_1,a'_1\ra \subseteq \Sigma$. 
    Finally, if such an $a$ does not exist then $\la a_1,a'_1\ra^\circ \subseteq \Sigma^c$, which ends the proof. 
\end{proof}

We further characterize when the cone $\la a_i,a'_i\ra$ is covered. 
It turns out that such a covering requires specific intersections with c-faces. 
\begin{proposition}\label{lem:a1ap1}
    Let $a_i,a'_i \in \bbr^n$. Then $\la a_i,a'_i\ra \subseteq \Sigma$ if and only if either $\dim (a_i,a'_i) \leq 1$ or $\dim (a_i,a'_i) = 2$ and there exists a c-face $F_i:=\la a_1,\dotsc,a_{i-1},a_{i+1},\dotsc,a_n\ra$ and a vector $p\neq 0$ 
    such that one of the following conditions occurs:  
    \[
        1.\ p \in \la a_i,a'_i\ra \cap F_i, \quad   
        2.\ p \in \la -a_i,a'_i\ra^\circ \cap F_i, \quad
        3.\ p \in \la a_i,-a'_i\ra^\circ \cap F_i \enspace .
    \]
    In all cases, one has $\la a'_i, a_1,\dotsc,a_n\ra \subseteq \Sigma$. 
\end{proposition}
\begin{proof}
    We fix $i$ to $1$ for clarity. 
    Suppose $\la a_1,a'_1\ra \subseteq \Sigma$ and that $\dim(a_1,a'_1)=2$.  
    Let $a \in \la a_1,a'_1\ra^\circ$, that is $a = \alpha_1 a_1 + \alpha'_1 a'_1$, $\alpha_1,\alpha'_1 > 0$. 
    Since $a \in \Sigma$, then either (i) $a \in \la a_1,\dotsc,a_n\ra$ or (ii) $a \in \la a'_1,\dotsc,a_n\ra$ for some c-face $F_1 = \la a_2,\dotsc,a_n\ra$. 
    For (i), one gets $\alpha_1 a_1 + \alpha'_1 a'_1 = \sum_{i=1}^n \beta_i a_i$ with $\beta_1,\dotsc,\beta_n \geq 0$. Thus $(\alpha_1 - \beta_1) a_1 + \alpha'_1 a'_1 = \sum_{i=2}^n \beta_i a_i$. 
    Let $p=(\alpha_1 - \beta_1) a_1 + \alpha'_1 a'_1$. Then $p\neq 0$ as otherwise $\dim(a_1,a'_1) \leq 1$. 
    \begin{itemize}
        \item If $\alpha_1 - \beta_1 \geq 0$ then $p\in\la a_1,a'_1\ra \cap F_1$ (case 1 of the statement). 
        \item If $\alpha_1 - \beta_1 < 0$ then 
    $p \in \la -a_1,a'_1\ra^\circ \cap F_1$ (case 2 of the statement).
    \end{itemize}
    The same discussion holds by swapping $a_1$ and $a'_1$ for (ii), leading to case 3. 
    For the converse, if $\dim(a_1,a'_1)\leq 1$ then it is immediate that $\la a_1,a'_1\ra \subseteq \Sigma$. Otherwise, one of the stated conditions holds. If the first condition holds, then by Lemma~\ref{lem:front} $\la a_1,a'_1,a_2,\dotsc,a_n\ra \subseteq \la a_1,a_2,\dotsc,a_n\ra \cup \la a'_1,a_2,\dotsc,a_n\ra \subseteq \Sigma$. If the second condition holds, then by Lemma~\ref{lem:side}, $\la a_1,a'_1,a_2,\dotsc,a_n\ra \subseteq \la a_1,\dotsc,a_n\ra \subseteq \Sigma$. 
    If the third condition holds, then by Lemma~\ref{lem:side} $\la a_1,a'_1,a_2,\dotsc,a_n\ra \subseteq \la a'_1,\dotsc,a_n\ra \subseteq \Sigma$ and $\la a_1,a'_1\ra \subseteq \Sigma$, as stated. 
\end{proof}




 

\subsection{Surrounding}\label{sec:self}

When the cone $\la a_i,a'_i\ra$ is degenerate, its front, side and back intersections with c-faces become ill-defined. We argue that these intersections are only relevant whenever $\la a_i,a'_i\ra$ is non-degenerate. 
We have already seen in Lemma~\ref{lem:inseparable} that $a_i \not\simeq a'_i$ for all $i$ is a necessary condition for $\Sigma$ to be covering. We shall see next that when $\la a_i,a'_i\ra$ is flat, we can drop the pair altogether and reduce the dimension of the Q-covering problem by $1$. 

\begin{proposition}\label{prop:reduction}
Assume that $a_i \neq 0$ and suppose that $\la a_i,a'_i \ra$ is flat. 
Let $a_i^\perp$ denote the hyperplane orthogonal to $a_i$, and let $\bar{v}$ denote the orthogonal projection of a vector $v$ onto $a_i^\perp$. 
Then $\Sigma = \bbr^n$ if and only if the sum $\{\ba_1,\ba'_1\}\oplus \dotsb\oplus \{\ba_{i-1},\ba'_{i-1}\}\oplus\{\ba_{i+1},\ba'_{i+1}\}\oplus \dotsb \oplus \{\ba_n, \ba'_n\}$ covers $a_i^\perp$ (which is isomorphic to $\bbr^{n-1}$).  
\end{proposition}
\begin{proof}
For clarity we fix $i$ to $1$. 
Let $h \in a_1^\perp$. Since $\bbr^n \subseteq \Sigma$, the vector $x_h := \alpha_1 a_1 + h$ with $\alpha_1 \in \bbr$, belongs to a c-cone $\la a_1,\dotsc,a_n\ra$ say. Thus $x_h = \sum_{i=1}^n \lambda_i a_i$,  $\lambda_1,\dotsc,\lambda_n \geq 0$, and $h = x_h - \alpha_1 a_1 = (\lambda_1 - \alpha_1) a_1 + \sum_{i \geq 2} \lambda_i a_i$, and $h = \bar{h} = \sum_{i \geq 2} \lambda_i \ba_i$ as stated. The same occurs if $x_h \in \la a'_1,a_2,\dotsc,a_n\ra$ since $a_1 \simeq -a'_1$. 
To prove the converse, let $x\in \bbr^n$ and decompose $x = \alpha_1 a_1 + h$ with $h \in a_1^\perp$. 
Thus $h$ belongs to some cone $\la \ba_2,\dotsc,\ba_n\ra$. 
Equivalently, there exist $\lambda_i \geq 0$ such that $h = \sum_{i \geq 2} \lambda_i \ba_i$. 
We have $a_i = \gamma_i a_1 + \ba_i$, $\gamma_i \in \bbr$, for each $i$. Therefore 
\[
\begin{split}
x &= \alpha_1 a_1 + h = \alpha_1 a_1 + \sum_{i \geq 2} \lambda_i \ba_i 
  = \alpha_1 a_1 + \sum_{i \geq 2} \lambda_i (a_i - \gamma_i a_1) \\ 
  &= \Bigl(\underbrace{\alpha_1 - \sum_{i \geq 2} \lambda_i \gamma_i}_\alpha \Bigr) a_1 + \sum_{i \geq 2} \lambda_i a_i  
\end{split}
\]
If $\alpha \geq 0$ then $x \in \la a_1,\dotsc,a_n\ra \subseteq \Sigma$ as required. If $\alpha < 0$ then there exists $\alpha' \geq 0$ such that $\alpha a_1 = \alpha' a'_1$ (using $a_1 \simeq -a'_1$) and $x \in \la a'_1,\dotsc,a_n\ra$ proving that $x \in \Sigma$ as well. 
\end{proof}

Self surrounding (cf. Definition~\ref{def:self}) is in fact equivalent to a Q-covering problem in lower dimension. 
\begin{proposition}\label{prop:self}
    Assume that $a_i \neq 0$ and let $\bar{v}$ denote the orthogonal projection of a vector $v$ onto $a_i^\perp$, the hyperplane orthogonal to $a_i$. 
    Then $a_i$ is self-surrounded if and only if $\{\ba_1,\ba'_1\}\oplus \dotsb\oplus \{\ba_{i-1},\ba'_{i-1}\}\oplus\{\ba_{i+1},\ba'_{i+1}\}\oplus \dotsb \oplus \{\ba_n, \ba'_n\}$ defines a Q-covering of $a_i^\perp$. 
\end{proposition}
\begin{proof}
    Let's fix $i$ to $1$. 
    If there exists a neighborhood $U$ of $a_1$ such that $U \subseteq \Sigma(a_1)$, then its projection $\bar{U}$ onto $a_1^\perp$ contains $0$ in its relative interior and must be included by definition of $\Sigma(a_1)$ in $\{\ba_2,\ba'_2\} \oplus \dotsb \oplus \{\ba_n,\ba'_n\}$. 
    Any nonzero vector in $a_1^\perp$ has a representative (with respect to $\simeq$) in $\bar{U}$ and is therefore covered. Thus $a_1^\perp \subseteq \{\ba_2,\ba'_2\} \oplus \dotsb \oplus \{\ba_n,\ba'_n\}$. 
    For the converse, consider $\Sigma' = \{a_1,-a_1\}\oplus\{a_2,a'_2\} \oplus \dotsb \oplus \{a_n,a'_n\}$. By Proposition~\ref{prop:reduction}, $\bbr^n \subseteq \Sigma'$ and therefore $a_1$ is surrounded with respect to $\Sigma'$. 
    Moreover, $a_1$ can only be self-surrounded (with respect to $\Sigma'$) since if it belongs to any c-cone in $\operatorname{cones}(\Sigma'(-a_1))$ then such a c-cone must be degenerate and therefore it does not effectively contribute in surrounding $a_1$. 
    Thus there exists a neighborhood $U$ of $a_1$ such that $U \subset \Sigma'(a_1)$. 
    By definition of $\Sigma'$, $\Sigma'(a_1) = \Sigma(a_1)$. Thus $U \subset \Sigma(a_1)$ and $a_1$ is self-surrounded with respect to $\Sigma$ as desired. 
\end{proof}
\begin{remark}
    In light of Proposition~\ref{prop:self}, Proposition~\ref{prop:reduction} can be reformulated as follows. Suppose that there exists an index $i$ such that $a_i \neq 0$ and that $\la a_i,a'_i\ra$ is flat. Then $\bbr^n$ is covered if and only if $a_i$ is self surrounded. 
\end{remark}

In general, the Q-covering problem does not enjoy an inductive property, that is, if $\Sigma = \bbr^n$, it is not necessary that $a_i$ is self surrounded for all $i$.  
It is worth mentioning that, in~\cite{cottle1980completely}, a subclass of Q-matrices, coined \emph{completely Q-matrices}, was introduced in which the inductive nature is preserved. 

Similar to self surrounding, lazy surrounding (cf. Definition~\ref{def:lazy}) can be also seen as a Q-covering problem in dimension $n-1$. The next proposition is stated for $a_1$ for clarity. The statement holds however for any $a_i$ with appropriate changes. 

%
%


\begin{proposition}\label{lem:lazy} 
    Suppose $a_1$ belongs to the topological interior of a c-cone $C=\la a'_1,a_2,\dotsc,a_n \ra$. 
    Then $\{\ba'_1,\ba_2\}\oplus\dotsb\oplus\{\ba'_1,\ba_n\}$ defines a Q-covering of $a_1^\perp$. \\
    The converse holds only when  $\det(\begin{smallmatrix}
        a_1 & a_2 & \dots & a_n
    \end{smallmatrix}) \det(\begin{smallmatrix}
        a'_1 & a_2 & \dots & a_n
    \end{smallmatrix}) > 0$.~\footnote{In this case, the hyperplane $(a_2,\dotsc,a_n)$ is sometimes said to be \emph{reflective}.} 
\end{proposition}
\begin{proof}
    Suppose $a_1$ belongs to the topological interior of a c-cone $C=\la a'_1,a_2,\dotsc,a_n \ra$, 
    then $a_1 \neq 0$ and $0 \in \la -a_1,a'_1,a_2,\dotsc,a_n\ra^\circ$. By Lemma~\ref{lem:flatcone}, $\la -a_1,a'_1,a_2,\dotsc,a_n\ra = (-a_1,a'_1,a_2,\dotsc,a_n) = (a'_1,a_2,\dotsc,a_n) = \bbr^n$ (because $C$ is non-degenerate) and 
    \[
    \la -a_1,a'_1, a_2,\dotsc,a_n\ra = C \cup \la -a_1, a_2,\dotsc,a_n \ra \cup \cup_{i=2}^n\la -a_1, a'_1, a_2,\dotsc,a_{i-1},a_{i+1}\dotsc,a_n \ra = \bbr^n \ .
    \]
    Moreover, by Lemma~\ref{lem:partition}, $C = \la a_1,\dotsc,a_n \ra \cup \cup_{i=2}^n\la a'_1,a_1,a_2,\dotsc,a_{i-1},a_{i+1}\dotsc,a_n\ra$. 
    Let $\Sigma'=\{a_1,-a_1\}\oplus\{a'_1,a_2\}\oplus\dotsb\oplus\{a'_1,a_n\}$. Then by definition of $\Sigma'$ one has 
    \begin{align*}
        \Sigma'(a_1) &= \la a_1,\dotsc,a_n \ra \cup \cup_{i=2}^n\la a_1,a'_1,a_2,\dotsc,a_{i-1},a_{i+1}\dotsc,a_n\ra = C \\
         \Sigma'(-a_1) &= \la -a_1, a_2,\dotsc,a_n \ra \cup \cup_{i=2}^n\la -a_1, a'_1, a_2,\dotsc,a_{i-1},a_{i+1}\dotsc,a_n \ra 
    \end{align*}
    Therefore, $\Sigma'(a_1) \cup \Sigma'(-a_1)=\bbr^n$. By Proposition~\ref{prop:reduction}, $\{\ba'_1,\ba_2\}\oplus\dotsb\oplus\{\ba'_1,\ba_n\}$ defines a Q-covering of $a_1^\perp$.  
    Conversely, if $\{\ba'_1,\ba_2\}\oplus\dotsb\oplus\{\ba'_1,\ba_n\}$ defines a Q-covering of $a_1^\perp$ then either $a_1$ or $-a_1$ is in the topological interior of $C$. 
    To ensure the former, $a_1$ and $a'_1$ must not be separated by the hyperplane $(a_2,\dotsc,a_n)$ which is equivalent to saying that  $\det(\begin{smallmatrix}
        a_1 & a_2 & \dots & a_n
    \end{smallmatrix})$ and $\det(\begin{smallmatrix}
        a'_1 & a_2 & \dots & a_n
    \end{smallmatrix})$ are both nonzero and have the same sign (if $a_1$ and $a'_1$ are separated by the hyperplane $(a_2,\dotsc,a_n)$ then surely $a_1 \not\in \la a'_1,a_2,\dotsc,a_n\ra$). 
\end{proof}

Both Propositions~\ref{prop:self} and~\ref{lem:lazy} will be used in Section~\ref{sec:surrounding} to arrive at an algebraic characterization of Q-matrices in dimension $3$.  
\section{Q-covering for  \texorpdfstring{$n=3$}{n=3}}
\label{sec:n=3}
A cone in $\operatorname{cones}(\Gamma) \setminus \operatorname{cones}(\Sigma)$ must be necessarily rooted at both $a_i$ and $a'_i$ for some indices $i$. Following~\cite{Garcia1983}, when such an index is unique, the cone is called an \emph{almost c-cone}. 
Assuming $\Gamma=\bbr^n$, it is clear that $\Sigma$ is covering if and only if all minimal cones are covered (necessity is obvious and sufficiency is a corollary of Proposition~\ref{prop:cara1}). 
This section investigates in depth the case $n=3$ with the objective of a complete understanding of when minimal cones are not covered. For this dimension, all cones in $\operatorname{cones}(\Gamma) \setminus \operatorname{cones}(\Sigma)$ are almost c-cones and have the form $\la a_i, a'_i, a_j\ra$, $j \neq i$. Thus for $n=3$, assuming $\Gamma=\bbr^3$, $\Sigma$ is covering if and only if all (full) almost c-cones are covered. 
This investigation leads to a complete understating of holes for this dimension (Theorem~\ref{thm:Gcovered}) as well as a characterization of the Q-covering problem for $n=3$ (Theorem~\ref{thm:symbolic}) which is amenable to an algebraic characterization as we shall see in Section~\ref{sec:surrounding}. 
It is worth mentioning that Corollary~\ref{coro:char} reduces the Q-covering problem to the surrounding of the vectors in $\Sigma_0$. It strengthens~\cite[Theorem 4.7]{Garcia1983} by dropping the \emph{strong non-degeneracy} assumption (which requires the non-degeneracy of all c-cones and almost c-cones). 
We start by a useful sufficient condition for minimal cones to be covered. 



\begin{proposition}\label{prop:ak}
    Let $i,j,k$ denote three distinct indices and let $G = \la a_i,a'_i,a_j\ra \subset \bbr^3$ denote a minimal cone. 
    Assume $\la a_i, a'_i\ra$ is covered. 
    If $-a_k \not\in G$ then $G \subseteq \Sigma$. 
\end{proposition}
\begin{proof}
    The proof is by case distinction following a partition of $\bbr^3$ suggested by $G$. Namely, $-G=\la -a_i,-a'_i,-a_j\ra$, the interior of $7$ full cones $C_i$, the relative interior of $9$ facets $F_i$ and $3$ rays $R_i$ where $C_1=G^\circ$, $C_2=\la a'_i,-a_i,a_j\ra^\circ$, $C_3=\la -a'_i,a_i,a_j\ra^\circ$, $C_4=\la -a_i,-a'_i,a_j\ra^\circ$, $C_5=\la a_i,a'_i,-a_j\ra^\circ$, $C_6=\la a'_i,-a_i,-a_j\ra^\circ$, and $C_7=\la -a'_i,a_i,-a_j\ra^\circ$; $F_1=\la a_i, a_j\ra^\circ$, $F_2=\la a_i, a'_i\ra^\circ$,  $F_3=\la a'_i, a_j\ra^\circ$, $F_4=\la -a_i, a'_i\ra^\circ$, $F_5=\la -a_i, a_j\ra^\circ$, $F_6=\la -a'_i, a_j\ra^\circ$, $F_7=\la a_i, -a'_i\ra^\circ$, $F_8=\la a_i, -a_j\ra^\circ$, and $F_9=\la a'_i, -a_j\ra^\circ$; $R_1=\la a_i\ra$, $R_2=\la a'_i\ra$, and $R_3=\la a_j\ra$. Cf. Fig.~\ref{fig:decomposition}.  
    By minimality of $G$, $a_k$ cannot be in $C_1,F_1,F_2,F_3$. 
    If $a_k \in C_5,C_6,C_7,F_4,F_7,F_8,F_9,R_1,R_2$, then $G \subseteq \la a_i,a'_i,a_j,a_k\ra \subseteq \Sigma$ where the last inclusion is a consequence of Proposition~\ref{lem:a1ap1}. 
    If $a_k \in C_4,F_5,F_6,R_3$, then $G \subseteq G'=\la a_i,a'_i,a_k\ra$. Since $\la a_i,a'_i\ra$ is covered, by Proposition~\ref{lem:a1ap1}, a c-face must intersect $\la a_i,a'_i \ra$ or $\la -a_i,a'_i \ra^\circ$ or $\la a_i,-a'_i \ra^\circ$ (such a c-face must have either $a_j$ or $a_k$ as generator) and either $G \subseteq \Sigma$ or $G' \subseteq \Sigma$. In both cases, $G$ is covered. 
    If $a_k \in C_2$ then $G \subseteq \la a_i,a_j,a_k\ra \cup G'$. Since $\la a_i,a'_i\ra$ is covered, as discussed above, a c-face having either $a_j$ or $a_k$ as generator satisfies the hypothesis of Proposition~\ref{lem:a1ap1}, and either $G$ is covered or $G'$ is covered. Thus, in all cases $G \subseteq \Sigma$. The same holds when $a_k \in C_3$ by swapping $a_i$ and $a'_i$.     
\end{proof}

\begin{figure}[bt]
     \centering
\begin{tikzpicture}[scale=.7]
\draw[thick] (0,0,0)--(0,0,2.5) node[below] {$R_1=\la a_i\ra$}; 
\draw[thick] (0,0,0)--(2.5,0,0) node[below] {$R_2=\la a'_i\ra$};  
\draw[thick] (0,0,0)--(0,1.5,0) node[left] {$R_3=\la a_j\ra$}; 
\draw[-stealth] (0,0,0)--(-2,0,0); 
\draw[-stealth] (0,0,0)--(0,0,-2); 

\coordinate[] (ai) at (0,0,2);
\coordinate[] (mai) at (0,0,-2);
\coordinate[] (api) at (2,0,0);
\coordinate[] (mapi) at (-2,0,0);
\coordinate[] (aj) at (0,1,0);

\draw[] (1.2,0,1.2) node[below] {$F_2$}; 
\coordinate[] (F2O) at (0.2,0,0.2); 
\coordinate[] (F2ai) at (2+.2,0,0.2);
\coordinate[] (F2api) at (0.2,0,2+.2); 

\coordinate[label=$F_4$] (F4) at (1.2,0,-1.2);
\coordinate[] (F4O) at (0.2,0,-0.2);
\coordinate[] (F4ai) at (2+0.2,0,-0.2);
\coordinate[] (F4mapi) at (0.2,0,-2-0.2);

\coordinate[label=$F_5$] (F5) at (0,0.7,-1.2);
\coordinate[] (F5O) at (0,0.2,-0.2);
\coordinate[] (F5mai) at (0,0.2,-2-0.2);
\coordinate[] (F5aj) at (0,1+.2,-0.2);

\draw[] (-1.2,0.7,0) node[above] {$F_6$};
\coordinate[] (F6O) at (-0.2,0.2,0);
\coordinate[] (F6mapi) at (-2-.2,0.2,0);
\coordinate[] (F6aj) at (-0.2,1+0.2,0);

\draw[] (-1.2,0,1.2) node[below] {$F_7$}; 
\coordinate[] (F7O) at (-0.2,0,0.2); 
\coordinate[] (F7api) at (-0.2,0,2+0.2);
\coordinate[] (F7mai) at (-2-0.2,0,0.2);

\draw[fill=gray,opacity=0.3] (F2O)--(F2ai)--(F2api)--(F2O);
\draw[fill=gray,opacity=0.3] (F4O)--(F4ai)--(F4mapi)--(F4O);
\draw[fill=gray,opacity=0.3] (F5O)--(F5mai)--(F5aj)--(F5O);
\draw[fill=gray,opacity=0.3] (F6O)--(F6mapi)--(F6aj)--(F6O);
\draw[fill=gray,opacity=0.3] (F7O)--(F7api)--(F7mai)--(F7O);
\end{tikzpicture}
         \qquad \qquad
\begin{tikzpicture}[scale=.7]
\draw[-stealth] (0,0,0)--(0,0,2) node[left] {$a_i$}; 
\draw[-stealth] (0,0,0)--(2,0,0) node[right] {$a'_i$};  
\draw[-stealth] (0,0,0)--(0,1.5,0) node[left] {$a_j$}; 
\draw[-stealth] (0,0,0)--(-2,0,0); 
\draw[-stealth] (0,0,0)--(0,0,-2); 
\draw[-stealth] (0,0,0)--(0,-1.5,0); 

\coordinate[] (ai) at (0,0,2);
\coordinate[] (mai) at (0,0,-2);
\coordinate[] (api) at (2,0,0);
\coordinate[] (mapi) at (-2,0,0);
\coordinate[] (aj) at (0,1,0);
\coordinate[] (maj) at (0,-1,0);

\draw[] (0,0.7,1.2) node[left] {$F_1$}; 
\coordinate[] (F1O) at (0,0.2,0.2); 
\coordinate[] (F1ai) at (0,0.2,2+.2);
\coordinate[] (F1aj) at (0,1+.2,0.2); 

\draw[] (1.2,0.7,0) node[above] {$F_3$}; 
\coordinate[] (F3O) at (0.2,0.2,0); 
\coordinate[] (F3api) at (2+.2,0.2,0);
\coordinate[] (F3aj) at (0.2,1+.2,0); 

\draw[] (0,-.7,1.2) node[below] {$F_8$}; 
\coordinate[] (F8O) at (0,-0.2,0.2); 
\coordinate[] (F8ai) at (0,-0.2,2+.2);
\coordinate[] (F8maj) at (0,-1-.2,0.2); 

\draw[] (1.2,-0.7,0) node[below] {$F_9$}; 
\coordinate[] (F9O) at (0.2,-0.2,0); 
\coordinate[] (F9api) at (2+.2,-0.2,0);
\coordinate[] (F9maj) at (0.2,-1-.2,0); 

\draw[fill=gray,opacity=0.3] (F1O)--(F1ai)--(F1aj)--(F1O);
\draw[fill=gray,opacity=0.3] (F3O)--(F3api)--(F3aj)--(F3O);
\draw[fill=gray,opacity=0.3] (F8O)--(F8ai)--(F8maj)--(F8O);
\draw[fill=gray,opacity=0.3] (F9O)--(F9api)--(F9maj)--(F9O);

\end{tikzpicture}
         \qquad \qquad
\begin{tikzpicture}[scale=.7]
\draw[-stealth] (0,0,0)--(0,0,2) node[left] {$a_i$}; 
\draw[-stealth] (0,0,0)--(2,0,0) node[right] {$a'_i$};  
\draw[-stealth] (0,0,0)--(0,1.5,0) node[right] {$a_j$}; 
\draw[-stealth] (0,0,0)--(-2,0,0); 
\draw[-stealth] (0,0,0)--(0,0,-2); 
\draw[-stealth] (0,0,0)--(0,-1.5,0); 

\coordinate[] (ai) at (0,0,2);
\coordinate[] (mai) at (0,0,-2);
\coordinate[] (api) at (2,0,0);
\coordinate[] (mapi) at (-2,0,0);
\coordinate[] (aj) at (0,1,0);
\coordinate[] (maj) at (0,-1,0);

\draw[] (2.6/3,1.6/3,2.6/3) node[] {$C_1$}; 
\coordinate[] (C1O) at (0.2,0.2,0.2); 
\coordinate[] (C1ai) at (0.2,0.2,2+.2);
\coordinate[] (C1api) at (2+.2,0.2,0.2);
\coordinate[] (C1aj) at (0.2,1.2,0.2); 

\draw[] (2.6/3,1.6/3,-2.6/3) node[above] {$C_2$};
\coordinate[] (C2O) at (0.2,0.2,-0.2); 
\coordinate[] (C2mai) at (0.2,0.2,-2-.2);
\coordinate[] (C2api) at (2+.2,0.2,-0.2);
\coordinate[] (C2aj) at (0.2,1.2,-0.2); 

\draw[] (-2.6/3,1.6/3,2.6/3) node[] {$C_3$}; 
\coordinate[] (C3O) at (-0.2,0.2,0.2); 
\coordinate[] (C3ai) at (-0.2,0.2,2+.2);
\coordinate[] (C3mapi) at (-2-.2,0.2,0.2);
\coordinate[] (C3aj) at (-0.2,1.2,0.2); 

\draw[] (-2.6/3,1.6/3,-2.6/3) node[above] {$C_4$}; 
\coordinate[] (C4O) at (-0.2,0.2,-0.2); 
\coordinate[] (C4mai) at (-0.2,0.2,-2-.2);
\coordinate[] (C4mapi) at (-2-.2,0.2,-0.2);
\coordinate[] (C4aj) at (-0.2,1.2,-0.2); 

\draw[] (2.6/3,-1.6/3,2.6/3) node[below] {$C_5$}; 
\coordinate[] (C5O) at (0.2,-0.2,0.2); 
\coordinate[] (C5ai) at (0.2,-0.2,2+.2);
\coordinate[] (C5api) at (2+.2,-0.2,0.2);
\coordinate[] (C5maj) at (0.2,-1.2,0.2); 

\draw[] (2.6/3,-1.6/3,-2.6/3) node[below right] {$C_6$};
\coordinate[] (C6O) at (0.2,-0.2,-0.2); 
\coordinate[] (C6mai) at (0.2,-0.2,-2-.2);
\coordinate[] (C6api) at (2+.2,-0.2,-0.2);
\coordinate[] (C6maj) at (0.2,-1.2,-0.2); 

\draw[] (-2.6/3,-1.6/3,2.6/3) node[below] {$C_7$}; 
\coordinate[] (C7O) at (-0.2,-0.2,0.2); 
\coordinate[] (C7ai) at (-0.2,-0.2,2+.2);
\coordinate[] (C7mapi) at (-2-.2,-0.2,0.2);
\coordinate[] (C7maj) at (-0.2,-1.2,0.2); 

\draw[fill=gray,opacity=0.3] (C1ai)--(C1api)--(C1aj)--(C1ai);
\draw[fill=gray,opacity=0.3] (C2mai)--(C2api)--(C2aj)--(C2mai);
\draw[fill=gray,opacity=0.3] (C3ai)--(C3mapi)--(C3aj)--(C3ai);
\draw[fill=gray,opacity=0.3] (C4mai)--(C4mapi)--(C4aj)--(C4mai);
\draw[fill=gray,opacity=0.3] (C5ai)--(C5api)--(C5maj)--(C5ai);
\draw[fill=gray,opacity=0.3] (C6mai)--(C6api)--(C6maj)--(C6mai);
\draw[fill=gray,opacity=0.3] (C7ai)--(C7mapi)--(C7maj)--(C7ai);
\end{tikzpicture}
        \caption{Decomposition of $\bbr^3$ in the proof of Proposition~\ref{prop:ak}.}
        \label{fig:decomposition}
\end{figure}

The following proposition is akin to the case of the so called self intersecting 3 \emph{starlike components} \cite[Fig. 1]{KELLY1979175}.  
We give below a direct (analytic) proof in line with our purposes. 
\begin{proposition}\label{prop:starlike2}
Let $G=\la a_i,a'_i,a_j\ra$ denote a full cone in $\bbr^3$ and 
suppose that there exists a nonzero vector $q \not\simeq a_j$ in $\la a_i,a'_j\ra \cap \la a'_i,a_j\ra$. Let $Q = \la a_i,a_j,q\ra$. 
If $-a_k \in Q$ then $Q^\circ \subseteq \Sigma(a_k)^c$. 
\end{proposition}
\begin{proof}
    For clarity we fix $\{i,j,k\}$ to $\{1,2,3\}$. It suffices to permute the indices to match any other configuration. 
    We first observe that $Q$ is non-degenerate because $q\neq 0$, $q \in \la a'_1,a_2\ra$ and $q \not\simeq a_2$. 
    By hypothesis, (i) $-a_3 = \lambda_1 a_1 + \lambda_2 a_2 + \lambda q$ with $\lambda_1,\lambda_2,\lambda \geq 0$. 
    Moreover, there exist $\gamma_1,\gamma'_1\geq 0$ and $\gamma_2,\gamma'_2 > 0$ such that (ii) $q = \gamma'_1 a'_1 + \gamma_2 a_2$ and (iii) $q = \gamma_1 a_1 + \gamma'_2 a'_2$. 
    Using (i) and (ii), one gets $-a_3 = \lambda_1 a_1 + \lambda_2 a_2 + \lambda (\gamma'_1 a'_1 + \gamma_2 a_2)$ or equivalently $(\lambda_2 + \lambda \gamma_2) a_2 + a_3 = \lambda_1 (-a_1)  + \lambda \gamma'_1 (-a'_1)$. Thus $\la -a_1,-a'_1\ra \cap \la a_2,a_3\ra$ is not reduced to zero. By Lemma~\ref{lem:back}, $G^\circ \subseteq (\la a_1,a_2,a_3\ra \cup \la a'_1,a_2,a_3 \ra)^c$ and $Q^\circ \subseteq \Sigma(a_2,a_3)^c$. 
    Moreover, let $x \in \la a_1,a_2,q\ra^\circ$, that is $x=\alpha_1 a_1 + \alpha_2 a_2 + \alpha q$, with $\alpha_1,\alpha_2,\alpha > 0$. 
    \begin{itemize}
    \item If $x \in \la a_1,q,a_3\ra$ then $x=\beta_1 a_1 + \beta q + \beta_3 a_3=\alpha_1 a_1 + \alpha_2 a_2 + \alpha q$ with $\beta_1,\beta,\beta_3 \geq 0$. 
    Thus $(\beta - \alpha) q + \beta_3 a_3 = (\alpha_1 - \beta_1) a_1 + \alpha_2 a_2$. 
    If $\beta_3 =0$ then $\alpha_2=0$ because $a_1,a_2,q$ are independent, contradicting $\alpha_2 > 0$. 
    Thus  $\beta_3 > 0$ and 
    (iv) $- \beta_3 a_3 = - (\alpha_1 - \beta_1) a_1 - \alpha_2 a_2 + (\beta - \alpha) q = \beta_3\lambda_1 a_1 + \beta_3\lambda_2 a_2 + \beta_3\lambda q$. Using (i) and (iv) one gets $-\alpha_2 = \beta_3\lambda_2$, a contradiction.  
    \item If $x \in \la a'_2,q,a_3\ra$ then $x=\beta'_2 a'_2 + \beta q + \beta_3 a_3=\alpha_1 a_1 + \alpha_2 a_2 + \alpha q$ with $\beta'_2,\beta,\beta_3 \geq 0$. 
    Thus $(\beta - \alpha) q + \beta_3 a_3 = \alpha_1 a_1 + \alpha_2 a_2 - \beta'_2 a'_2$. 
    If $\beta_3 =0$. Then $\beta - \alpha$ cannot be zero since $a_1,a'_2,a_2$ are independent ($q \not\simeq a_2$). Moreover, if $\beta - \alpha \neq 0$, $q \in \la a_1,a'_2\ra$ forces $\alpha_2$ to be zero contradicting $\alpha_2 > 0$. So $\beta_3 > 0$ and 
    (v) $- \beta_3 a_3 = - \alpha_1 a_1 - \alpha_2 a_2 + \beta'_2 a'_2 + (\beta - \alpha) q = \beta_3\lambda_1 a_1 + \beta_3\lambda_2 a_2 + \beta_3\lambda q$. 
    Using (iii), one substitutes $a'_2$ for a linear combination of $a_2$ and $q$. 
    Thus (i) and (v) lead to $-\alpha_1 = \beta_3\lambda_1$, a contradiction. 
    \end{itemize}
    Since $q \in \la a'_1,a_2\ra$, $\la a_1,a'_2,a_3\ra = \la a_1,q,a_3\ra \cup \la a'_2,q,a_3\ra$. Therefore $x \not\in \la a_1,a'_2,a_3\ra$. 
    Finally, if $x \in \la a'_1,a'_2,a_3\ra$ then $x=\beta'_1 a'_1 + \beta'_2 a'_2 + \beta_3 a_3=\alpha_1 a_1 + \alpha_2 a_2 + \alpha q$ with $\beta'_1,\beta'_2,\beta_3 \geq 0$. 
    Using $q = \gamma_1 a_1 + \gamma'_2 a'_2 = \gamma'_1 a'_1 + \gamma_2 a_2$, if $\gamma'_1 = 0$ then $q \simeq a_2$, contradicting the hypothesis on $q$. Thus $\gamma'_1 > 0$. If $\gamma'_2 = 0$ then $G$ would be degenerate, also a contradiction. Thus $\gamma'_2 > 0$. 
    One then gets (v) $\alpha_1 a_1 + \alpha_2 a_2 + \alpha q = \frac{\beta'_1}{\gamma'_1} (q-\gamma_2 a_2) + \frac{\beta'_2}{\gamma'_2} (q-\gamma_1 a_1) + \beta_3 a_3$. 
    If $\beta_3=0$ then $\alpha_1 + \tfrac{\beta'_2}{\gamma'_2}\gamma_1 = 0$ contradicting $\alpha_1 > 0$. Thus $\beta_3 > 0$ and using (i) and (v) one gets $-(\alpha_1 + \tfrac{\beta'_2}{\gamma'_2}\gamma_1) = \beta_3 \lambda_1$ which is impossible. Thus $x\not\in \la a'_1,a'_2,a_3\ra$ and $Q^\circ \subseteq \Sigma(a'_2,a_3)^c$. 
    Therefore $Q^\circ \subseteq \Sigma(a_2,a_3)^c \cap \Sigma(a'_2,a_3)^c = \Sigma(a_3)^c$. 
\end{proof}

Proposition~\ref{prop:a1a1pcov} showed that the cone $\la a_i, a'_i \ra$ enjoys a special property: it cannot be partially covered. 
When it is covered, Theorem~\ref{thm:Gcovered} below gives a necessary condition for a minimal cone $\la a_i,a'_i,a_j\ra$ to contain a hole. Its proof requires the following technical lemma. 

\begin{lemma}\label{lem:sideintersection2}
    Let $i,j,k$ denote distinct indices in $\{1,2,3\}$, $G = \la a_i,a'_i,a_j\ra$ denote a full cone in $\bbr^3$ such that $-a_k \in G$. 
    Suppose there exist vectors $p,q\neq 0$ such that $p\in \la -a_i, a'_i \ra^\circ \cap \la a'_j, a_k\ra^\circ$ and $q \in \la a_i,a'_j\ra^\circ \cap \la a'_i,a_j\ra^\circ$. 
    Then $-a_k \in \la a_i,a_j,q\ra$. 
\end{lemma}
\begin{proof}
    For clarity we fix $\{i,j,k\}$ to $\{1,2,3\}$. It suffices to permute the indices to match any other configuration. 
    Since $p\in \la -a_1, a'_1 \ra^\circ \cap \la a'_2, a_3\ra^\circ$, there exist $\lambda_1,\lambda'_1 > 0$ and $\lambda'_2,\lambda_3 > 0$ such that (i) $p=\lambda_1 (-a_1) + \lambda'_1 a'_1 = \lambda'_2 a'_2 + \lambda_3 a_3$. 
    Moreover since $q \in \la a_1,a'_2\ra^\circ \cap \la a'_1,a_2\ra^\circ$, there exist $\gamma'_1,\gamma_2,\gamma_1,\gamma'_2 > 0$, such that (ii)   
    $q = \gamma'_1 a'_1 + \gamma_2 a_2 = \gamma_1 a_1 + \gamma'_2 a'_2$. 
    We can eliminate $a'_2$ from (i) and (ii) to get 
    %
    \[
    -\gamma'_2\lambda_3 a_3 = \underbrace{(-\lambda'_2\gamma_1 + \gamma'_2\lambda_1)}_{\theta_1} a_1 + \underbrace{(-\gamma'_2\lambda'_1 + \lambda'_2\gamma'_1)}_{\theta'_1} a'_1 + \lambda'_2\gamma_2 a_2 
    \]
    with $\theta_1, \theta'_1 \geq 0$ (because $-a_3 \in G$). 
    Thus $\tfrac{\gamma'_1}{\lambda'_1} \geq \tfrac{\gamma'_2}{\lambda'_2} \geq \tfrac{\gamma_1}{\lambda_1}$ and $\rho_1 = \gamma'_1\lambda_1 - \lambda'_1\gamma_1 \geq 0$. 
    We further eliminate $a'_1$ from (i) and (ii) to get 
    \[\rho_1 a_1 + \lambda'_1\gamma_2 a_2 + \theta'_1 a'_2 = -\gamma'_1\lambda_3 a_3\]
    Thus $-a_3 \in G \cap \la a_1,a_2,a'_2\ra = \la a_1,a_2,q\ra$. 
\end{proof}

\begin{theorem}\label{thm:Gcovered}
    Let $G=\la a_i,a'_i,a_j\ra$, $j \neq i$, denote a minimal cone such that $\la a_i,a'_i\ra \subseteq \Sigma$. 
    If $G$ contains a hole $K$ then, up to swapping $a_i,a'_i$, there exists a nonzero vector $q \in \la a_i,a'_j\ra \cap \la a'_i,a_j\ra$ such that $K=\la a_i,q,a_j\ra^\circ$. 
\end{theorem} 
\begin{proof}
    If $-a_k$ or $-a'_k$ is not in $G$ then, by Proposition~\ref{prop:ak}, $G \subseteq \Sigma$. 
    By hypothesis $G$ contains a hole so we can assume that $-a_k,-a'_k \in G$. Since $\la a_i,a'_i\ra$ is covered, by Proposition~\ref{lem:a1ap1}, one c-face must intersect $\la a_i,a'_i\ra$ or $\la -a_i,a'_i\ra^\circ$ or $\la a_i,-a'_i\ra^\circ$. 
    In this case, such a c-face must have $a'_j$ as a generator (as both $\la a_j,a_k\ra$ and $\la a_j,a'_k\ra$ intersect $\la -a_i,-a'_i\ra$). 
    With respect to the partition of the space used in the proof of Proposition~\ref{prop:ak} (cf. Fig~\ref{fig:decomposition}), and since $G$ is a minimal cone, for $\la a_i,a'_i\ra$ to be covered, $a'_j$ must belong to one of the following $C_2,C_3,F_4,F_7,R_1,R_2$. 
    In all cases there exists a nonzero vector $q \not\simeq a'_j$ in $\la a_i,a'_j \ra \cap \la a'_i,a_j \ra$ (when $a'_j$ is in $C_2,F_4,R_2$) or in $\la a_i,a_j \ra \cap \la a'_i,a'_j \ra$ (when $a'_j$ is in $C_3,F_7,R_1$). 
    (Observe that both cases are symmetric by swapping $a_i$ and $a'_i$.)
    Suppose the former, and let $Q = \la a_i,a_j,q\ra$. 
    If $a'_j \in R_2$ or $a'_j \in F_4$ then $q \simeq a'_i$ and $Q=G$. 
    By Proposition~\ref{prop:starlike2}, $G^\circ = Q^\circ \subseteq \Sigma(a_k)^c \cap \Sigma(a'_k)^c = \Sigma^c$. 
    Suppose $a'_j \in C_2$, then there exists $p \in \la -a_i,a'_i\ra^\circ \cap \la a'_j,a_k\ra^\circ$. By Lemma~\ref{lem:sideintersection2}, $-a_k \in Q$ and by Proposition~\ref{prop:starlike2}, $Q^\circ \subseteq \Sigma(a_k)^c$. 
    If $\la a'_j,a'_k \ra$ does not intersect $\la -a_i,a'_i\ra^\circ$ then it intersects $\la -a_i,-a'_i\ra$ and therefore $G^\circ \subseteq \Sigma(a'_k)^c$ by Lemma~\ref{lem:back}. 
    Thus $Q^\circ \subseteq \Sigma^c$. 
    Otherwise $\la a'_j,a'_k \ra$ intersect $\la -a_i,a'_i\ra^\circ$ and $-a'_k \in Q$ and by Proposition~\ref{prop:starlike2}, $Q^\circ \subseteq \Sigma(a'_k)^c$, $Q^\circ \subseteq \Sigma^c$. 
    As $Q^\circ$ is the maximal hole contained in $G$, it follows that $K=Q^\circ$ as stated. 
\end{proof}

\begin{corollary}\label{coro:char}
    $\Sigma = \bbr^3$ if and only if, for all $i$, both $a_i$ and $a'_i$ are surrounded.~\footnote{In \cite{KELLY1979175}, the authors relied heavily on visualization to characterize $3\times 3$ non-degenerate Q-matrices. Interested readers can find a proof in the same spirit in~\cite[Proposition 4]{CKThesis}.} 
\end{corollary}
\begin{proof}
    Necessity is immediate. 
    For sufficiency, we prove the contrapositive, that is if $\Sigma^c$ is non-empty then there exists an index $i$ for which either $a_i$ or $a'_i$ is not surrounded. Assume first that $\Gamma\subset \bbr^3$. 
    Then $\Gamma$ is a closed proper convex cone of $\bbr^3$ with a boundary that is non-empty. 
    If for all $a\in \Sigma_0$, $a\in \Gamma^{\circ}$ then $\Gamma \subseteq \Gamma^\circ \subseteq \Gamma$. So $\Gamma$ is both open and closed and its boundary must be empty, a contradiction.~\footnote{Another way to arrive at a contradiction would be to use the fact that $\bbr^n$ is a connected topological space, thus the only clopen subsets are $\bbr^n$ and its complement, excluding $\Gamma$.} Thus, there exists a vector $a \in \Sigma_0$ such that $a$ is a boundary ray of $\Gamma$ and therefore $a$ cannot be surrounded (if $U$ is a neighborhood of $a$, then $U \not\subseteq \Gamma$ and since $\Sigma \subseteq \Gamma$, $U \not\subseteq \Sigma$). 
    Next, assume that $\Gamma = \bbr^3$. 
    Since $\Sigma \subset \bbr^3$, by Proposition~\ref{prop:cara1}, there must exist a minimal cone $G = \la a_i,a'_i,a_j\ra$, $j \neq i$, such that $K=G \cap \Sigma^c$ is non-empty. 
    By Proposition~\ref{prop:a1a1pcov}, either $\la a_i, a'_i\ra^\circ \subseteq \Sigma^c$ or $\la a_i, a'_i\ra \subseteq \Sigma$. If the former holds then both $a_i$ and $a'_i$ are not surrounded. If the latter holds, then by Theorem~\ref{thm:Gcovered}, either $a_i$ is not surrounded or $a'_i$ is not surrounded.  
\end{proof}


Corollary~\ref{coro:char} strengthens~\cite[Theorem 4.7]{Garcia1983} by dropping the strong non-degenerate assumption. The latter result was moreover established using degree theory while in this work we solely used convex geometry. 
The statement of Corollary~\ref{coro:char} does not hold for $n>3$. In~\cite[Theorem 3]{morris1988counterexamples}, the author gives an example in $n=4$ where both $a_i$ and $a'_i$ are (lazily) surrounded for all $i$ without having a Q-covering. 
The (counter)example in $n=4$ thus shows that it is possible for a minimal cone to have surrounded generators while still having a hole in it. 

Another (more practical) issue with the statement of Corollary~\ref{coro:char} is that, in general, surrounding is not straightforward to transpose algebraically. 
Under the assumption $\Gamma = \bbr^3$, in the sequel we alleviate the need for checking surrounding: we show that self surrounding and lazy covering are enough to characterize the Q-covering for $n=3$. 
We start by proving some special cases before stating the main theorem (cf. Theorem~\ref{thm:symbolic} below). 
\begin{lemma}\label{lem:ap1a2}
    Suppose that $a_i \in \la a'_i,a_j\ra$, $i \neq j$. If $a'_i$ is self surrounded or lazily covered then all minimal cones rooted at $a_i,a'_i$ are covered. 
\end{lemma}
\begin{proof}
    For clarity we fix $(i,j)$ to $(1,2)$. (It suffices to permute the indices accordingly for any other configuration.) 
    Since $a_1 \in \la a'_1, a_2\ra$, then $\la a_1,a'_1,a_2\ra \subseteq \la a'_1,a_2\ra \subseteq \Sigma$, $\la a_1,a'_1,a_3\ra \subseteq \la a'_1,a_2,a_3\ra \subseteq \Sigma$, and $\la a_1,a'_1,a'_3\ra \subseteq \la a'_1,a_2,a'_3\ra \subseteq \Sigma$. 
    Only $G=\la a_1,a'_1,a'_2\ra$ is left out. 
    Assume $a'_1$ is self surrounded. 
    Thus $\la a_1,a'_1\ra$ is covered by Proposition \ref{prop:a1a1pcov}. Suppose $-a_3,-a'_3 \in G$, then by Proposition~\ref{prop:starlike2},  $G^\circ \subseteq \Sigma(a_3)^c \cap \Sigma(a'_3)^c = \Sigma^c$, contradicting the self surrounding of $a'_1$.  
    Thus either $-a_3 \not\in G$ or $-a'_3 \not\in G$ and by Proposition~\ref{prop:ak}, $G$ is covered as desired.
    Assume $a'_1 \in \Sigma(a_1)$. If $a'_1 \in \la a_1,a'_2,a_3\ra$, then $G \subseteq \la a_1,a'_2,a_3\ra \subseteq \Sigma$. The same occurs if $a'_1 \in \la a_1,a'_2,a'_3\ra$ by swapping $a_3$ and $a'_3$. 
    If $a'_1 \in \la a_1,a_2,a_3\ra$ then either $a_3 \simeq a'_1$ or $a_3 \in \la -a_1,a'_1\ra^\circ$. 
    The former makes $G$ a c-cone. If the latter occurs, then by Lemma~\ref{lem:side}, $G$ is covered. 
    The same discussion holds if $a'_1 \in \la a_1,a_2,a'_3\ra$ by swapping $a_3$ and $a'_3$. We thus proved that all minimal cones rooted at $a_1,a'_1$ are covered.      
\end{proof}

\begin{lemma}\label{lem:a2a3}
    Suppose $a_i \in \la a_j, a_k\ra^\circ$ with $i,j,k$ distinct. 
    Then all minimal cones rooted at $a_i,a'_i$ are covered. 
\end{lemma}
\begin{proof}
    For clarity we fix $i,j,k$ to $1,2,3$ respectively. 
    (It suffices to permute the indices accordingly for any other configuration.) 
    Since $a_1 \in \la a_2, a_3\ra$ then $\la a_1,a'_1,a_3\ra \subseteq \la a'_1,a_2,a_3\ra \subseteq \Sigma$ and $\la a_1,a'_1,a_2\ra \subseteq \la a'_1,a_2,a_3\ra \subseteq \Sigma$. 
    It remains to prove that $\la a_1,a'_1,a'_2\ra$ and $\la a_1,a'_1,a'_3\ra$ are covered. 
    We first prove that $G=\la a_1,a'_1,a'_2\ra$ is covered. 
    If $a_2 \simeq a'_2$, then $G$ is also covered as discussed above. 
    If $a_2 \simeq a'_1$ then $G \subseteq \la a'_1,a'_2,a_3\ra \subseteq \Sigma$. In addition $a_2 \not\simeq a_1$ (because $a_1 \in \la a_2,a_3\ra^\circ$). Using the minimality of $G$, we can thus assume in the sequel that $a_2 \not\in G$.  
    Suppose that $-a_3 \in G$, then there exists $\lambda_1,\lambda'_1,\lambda'_2 \geq 0$ such that $-a_3 = \lambda_1 a_1 + \lambda'_1 a'_1 + \lambda'_2 a'_2$. Moreover there exists $\gamma_2,\gamma_3 > 0$ such that $a_1 = \gamma_2 a_2 + \gamma_3 a_3$. By eliminating $a_3$ from both equations, one gets 
    \[
    \gamma_2 a_2 = (\gamma_3\lambda_1 + 1) a_1 + \gamma_3\lambda'_1 a'_1 + \gamma_3\lambda'_2 a'_2 
    \]
    and therefore $a_2 \in G$, a contradiction. Thus $-a_3 \not\in G$. 
    By hypothesis $\la a_1,a'_1\ra \subseteq \la a'_1,a_2,a_3\ra \subseteq \Sigma$. Thus Proposition~\ref{prop:ak} applies and $G \subseteq \Sigma$. 
    The exact same discussion holds to prove that $\la a_1,a'_1,a'_3\ra$ is covered by swapping the indices $2$ and $3$. 
\end{proof}

\begin{theorem}\label{thm:symbolic}
    Assume $\Gamma = \bbr^3$. Then 
    $\Sigma = \bbr^3$ if and only if, all vectors in $\Sigma_0$ are either self surrounded or lazily covered. (Observe that necessity holds for any finite dimension.) 
\end{theorem} 
\begin{proof} 
    (Necessity) If $\Sigma = \bbr^3$, then $a_i$ must be surrounded for all $i$. If $a_i$ is not self surrounded then necessarily $a_i \in \Sigma(a'_i)$ (otherwise $a_i$ cannot be surrounded). The same holds for $a'_i$. 
    (Sufficiency) We prove the contrapositive, i.e. if $\Sigma \subset \bbr^3$ then there exists a vector in $\Sigma_0$ which is not self surrounded nor lazily covered. 
    The proof is by contradiction. Suppose that $\Sigma \subset \bbr^3$ and all vectors in $\Sigma_0$ are either self surrounded or lazily covered. Thus, there exists a hole $K$. 
    Since $\Gamma = \bbr^3$, then by Proposition~\ref{prop:cara1} there exists a minimal cone $G=\la a_i,a'_i,a_j\ra$, $j \neq i$, such that $K=G \cap \Sigma^c$ is non-empty. 
    By Proposition~\ref{prop:a1a1pcov}, either $\la a_i,a'_i\ra^\circ \subset \Sigma^c$ or $\la a_i,a'_i\ra \subset \Sigma$. If the former holds, then $a_i$ (or $a'_i$) leads to a contradiction as it is not self surrounded nor lazily covered (cf. the proof of Proposition~\ref{prop:a1a1pcov}). If the latter holds, by Theorem~\ref{thm:Gcovered}, $a_i$, say, is not surrounded (otherwise it suffices to swap $a_i$ and $a'_i$ in what follows). 
    In particular $a_i$ is not self nor lazily surrounded. 
    Therefore $a_i \in \Sigma_2(a'_i)$ (i.e. a face of a c-cone rooted at $a'_i$). 
    If $a_i \in \la a'_i,a_j\ra$, $j \neq i$, then by Lemma~\ref{lem:ap1a2}, $G$ is covered. 
    If $a_i \in \la a_j,a_k\ra^\circ$, with $i,j,k$ distinct, then by Lemma~\ref{lem:a2a3}, $G$ is also covered. 
    Both cases lead to a contradiction since $K$ is non-empty. 
    Thus $K$ does not exist and $\Sigma=\bbr^3$. 
\end{proof}

Unlike Corollary~\ref{coro:char}, Theorem~\ref{thm:symbolic} is amenable to an \emph{algebraic characterization} of Q-matrices for $n=3$ as we shall detail in the next section. 
 
\section{Algebraic Characterization}\label{sec:surrounding}
Theorem~\ref{thm:symbolic} characterizes the Q-covering via three requirements: (1) Feasibility or S-matricity (i.e. $\Gamma = \bbr^3$), (2) self surrounding (i.e. the Q-covering of $a_i^\perp$), and (3) lazy covering (i.e. $a_i \in \Sigma(a'_i)$). 
We detail next how each requirement can be equivalently translated into sign conditions on the subdeterminants of the matrix $(\begin{smallmatrix}
    a_1 & a_2 & a_3 & a'_1 & a'_2 & a'_3
\end{smallmatrix})$. 

\subsection{S-matricity}\label{sec:smatsymbolic} 

Proposition~\ref{prop:gammacovering} provides an effective mean to characterize S-matrices for $n=3$. 
\begin{corollary}\label{coro:gamma} 
    Let $\Gamma = \la g_1,\dotsc,g_{6} \ra$. Then 
    $\Gamma = \bbr^3$ if and only if 
    \begin{itemize}
        \item either there are $4$ vectors such that $\bbr^3 = \la g_{i_1},\dotsc,g_{i_4}\ra$,  
        \item or there are $3$ vectors such that $\la g_{i_1},g_{i_2},g_{i_3}\ra$ is a plane that separates two other vectors $g_{i_4}$ and $g_{i_5}$, 
        \item or there are $2$ vectors such that $\la g_{i_1},g_{i_2}\ra$ is a line and the plane $g_{i_{1}}^\perp$ is equal to $\la \pi(g_{i_3}),\dotsc,\pi(g_{i_6})\ra$, where $\pi$ denotes the orthogonal projection onto the hyperplane $g_{i_1}^\perp$. 
    \end{itemize}
\end{corollary} 
\begin{proof}
 Sufficiency is immediate. For necessity, by Proposition~\ref{prop:gammacovering}, there exists $1 \leq m \leq 3$ such that $m+1$ vectors among $g_1,\dotsc,g_6$ span a flat of dimension $m$. The provided conditions enumerate all cases ($m=3$ first, then $m=2$, and finally $m=1$). 
\end{proof}

When the four cones $\la g_1,g_2,g_3 \ra$, $\la g_2,g_3,g_4 \ra$, $\la g_1,g_3,g_4 \ra$, and $\la g_1,g_2,g_4 \ra$ are full and have the same orientation, then $\Gamma = \la g_1,g_2,g_3,g_4\ra=\bbr^3$. The orientation can be retrieved using the sign of the determinant of the matrix formed by the generators. 
For instance, the orientation of the cone $\la g_1,g_2,g_3\ra$ is given by the sign of $\det(\begin{smallmatrix}
    g_1 & g_2 & g_3
\end{smallmatrix})$. 
However, while the cone, as a geometric object, is invariant under any permutation of its generators, the sign of the determinant is not. One thus has to be careful when ordering the vectors to get a \emph{coherent orientation} of the involved cones~\cite[Section 18]{invittopalg2014}. 
We do this by fixing a global order of the involved vectors and making sure that the orientation of the common facet of any two adjacent cones is inverted w.r.t. the fixed global order. 
For instance the cones spanned by the lists $\{g_1,g_2,g_3 \}$ and $\{ g_2,g_3,g_4 \}$ are adjacent cones having in common $\{ g_2,g_3\}$. 
With respect to the global ordering $g_1,g_2,g_3,g_4$, the lists $\{g_1,g_2,g_3\}$ and $\{g_3,g_2,g_4\}$ are coherently oriented.

Similarly, when $\bbr^3$ is nonnegatively spanned by five vectors, assuming any $4$ of them do not span the space, then three of them span a plane that separates the two remaining ones. This implies that the space is partitioned into six distinct cones that have the same orientation (i.e. all pairs of adjacent cones are coherently oriented). For a list of $5$ vectors, we misuse the $\oplus$ notation (cf. Def.~\ref{def:qmatrix}) and write $\Gamma$ as $\{g_1,g_2,g_3\} \oplus \{g_4,g_5\}$, where the first three vectors span the plane separating the two remaining ones.  
Finally, when the six vectors are required to span $\bbr^3$, they must form $3$ lines in a generic position (no one is in the plane formed by the two others). This implies that the space is partitioned into eight distinct cones that have the same orientation. In this case, we write $\Gamma$ as $\{g_1,g_2\} \oplus \{g_3,g_4\} \oplus \{g_5,g_6\}$ where each pair form a line. 
The notation here is suggestive of the Q-covering problem. In fact it is a very special case where one gets a partition of the space. Equivalently, $(\begin{smallmatrix}
    g_1 g_3 g_5
\end{smallmatrix})^{-1}(\begin{smallmatrix}
    g_2 g_4 g_6
\end{smallmatrix})$ is a P-matrix. 

We say that a vector is \emph{symbolic} if all its components are non-fixed reals or variables or \emph{symbols}. 
For instance the vector $(m_1,m_2,m_3)$ is symbolic whereas $(1,0,0)$ is not. 
Algorithm~\ref{alg:gamma4} makes explicit the conditions that the components of $4$ vectors have to satisfy to span $\bbr^3$. 
Algorithm~\ref{alg:gamma5} makes explicit the conditions that a list of $5$ vectors have to satisfy to span $\bbr^3$ assuming the first three vectors span a plane and the remaining two vectors are separated by that plane. To avoid checking whether the first three vectors actually form a plane, the algorithm relies on partitioning the space into $6$ coherently oriented cones. 
The so obtained conditions may thus have some redundancy with the ones obtained from Algorithm~\ref{alg:gamma4}. 
Finally, Algorithm~\ref{alg:gamma6} makes explicit the conditions that a list of $6$ vectors have to satisfy to span $\bbr^3$ (where all the vectors are required). As discussed earlier, we also implement the partition of the space into $8$ cones for simplicity at the cost of redundancy with the conditions provided by the two other algorithms. 

\begin{algorithm}[t]
  \DontPrintSemicolon
  \LinesNumbered
  \SetAlgoVlined
  \KwData{Four symbolic vectors $g_1,\dotsc,g_4$.} 
   $\{d_1,d_2,d_3,d_4\} \leftarrow \{\det(\begin{smallmatrix}
       g_1 & g_2 & g_3 
   \end{smallmatrix}), \det(\begin{smallmatrix}
       g_2 & g_1 & g_4
   \end{smallmatrix}), \det(\begin{smallmatrix}
       g_3 & g_4 & g_1
   \end{smallmatrix}), \det(\begin{smallmatrix}
       g_4 & g_3 & g_2
   \end{smallmatrix}) \}$ \; 
    \KwRet{$(\bigwedge_{i=1}^4 d_i > 0) \lor (\bigwedge_{i=1}^4 d_i < 0)$} 
  \caption{$\Gamma_4$: $\bbr^3 = \la g_1,g_2,g_3,g_4 \ra$ .}\label{alg:gamma4}
\end{algorithm}

\begin{algorithm}[t]
  \DontPrintSemicolon
  \LinesNumbered
  \SetAlgoVlined
  \KwData{Five symbolic vectors $g_1,\dotsc,g_5$.} 
   $\{d_1,d_2,d_3\} \leftarrow \{\det(\begin{smallmatrix}
       g_1 & g_2 & g_4 
   \end{smallmatrix}), \det(\begin{smallmatrix}
       g_2 & g_3 & g_4
   \end{smallmatrix}), \det(\begin{smallmatrix}
       g_3 & g_1 & g_4
   \end{smallmatrix}) \}$ \; 
   $\{d_4,d_5,d_6\} \leftarrow \{\det(\begin{smallmatrix}
       g_1 & g_3 & g_5
   \end{smallmatrix}), \det(\begin{smallmatrix}
       g_3 & g_2 & g_5
   \end{smallmatrix}), \det(\begin{smallmatrix}
       g_2 & g_1 & g_5
   \end{smallmatrix}) \}$ \; 
    \KwRet{$(\bigwedge_{i=1}^6 d_i > 0) \lor (\bigwedge_{i=1}^6 d_i < 0)$} 
  \caption{$\Gamma_5$: $\bbr^3 = \la g_1,\dotsc,g_5 \ra = \{g_1,g_2,g_3\}\oplus\{g_4,g_5\}$ .}\label{alg:gamma5}
\end{algorithm}

\begin{algorithm}[t]
  \DontPrintSemicolon
  \LinesNumbered
  \SetAlgoVlined
  \KwData{Six symbolic vectors $g_1,\dotsc,g_6$.} 
   $\{d_1,d_2,d_3,d_4\} \leftarrow \{\det(\begin{smallmatrix}
       g_1 & g_3 & g_5 
   \end{smallmatrix}), \det(\begin{smallmatrix}
       g_3 & g_2 & g_5
   \end{smallmatrix}), \det(\begin{smallmatrix}
       g_2 & g_4 & g_5
   \end{smallmatrix}), \det(\begin{smallmatrix}
       g_4 & g_1 & g_5
   \end{smallmatrix}) \}$ \; 
   $\{d_5,d_6,d_7,d_8\} \leftarrow \{\det(\begin{smallmatrix}
       g_1 & g_4 & g_6
   \end{smallmatrix}), \det(\begin{smallmatrix}
       g_4 & g_2 & g_6
   \end{smallmatrix}), \det(\begin{smallmatrix}
       g_2 & g_3 & g_6
   \end{smallmatrix}), \det(\begin{smallmatrix}
       g_3 & g_1 & g_6
  \end{smallmatrix}) \}$ \; 
    \KwRet{$(\bigwedge_{i=1}^8 d_i > 0) \lor (\bigwedge_{i=1}^8 d_i < 0)$} 
  \caption{$\Gamma_6$: $\bbr^3 = \la g_1,\dotsc,g_6 \ra = \{g_1,g_2\}\oplus\{g_3,g_4\}\oplus\{g_5,g_6\}$ .}\label{alg:gamma6}
\end{algorithm}

Given $6$ symbolic vectors $g_1,\dotsc,g_6$, one can characterize $\Gamma = \la g_1,\dotsc,g_6\ra = \bbr^3$ as follows: apply Algorithm~\ref{alg:gamma4} to any sublist of four vectors among the ones provided, apply Algorithm~\ref{alg:gamma5} to all distinct pairs of sublists of $3$ and $2$ vectors, and finally apply Algorithm~\ref{alg:gamma6} to all distinct tuples formed each by three pairs of vectors. 
 
\subsection{Self Surrounding}\label{sec:planar}

By Proposition~\ref{prop:self}, self surrounding for $n=3$ amounts to characterizing the Q-covering for $n=2$. 
The following theorem is the analogue of Corollary~\ref{coro:char} for dimension $2$. It relies solely on surrounding which is easy to characterize for this dimension, alleviating the need for feasibility.  
\begin{theorem}\label{thm:n2surrounding}
    $\Sigma = \bbr^2$ if and only if, for all $i$, $a_i$ and $a'_i$ are surrounded. 
\end{theorem}
\begin{proof}
    Necessity is immediate. For sufficiency, we prove the contrapositive. Suppose $\Sigma \subset \bbr^2$. If $\Gamma \subset \bbr^2$ then by~\cite[Corollary 18.3.1]{rockafellar1997convex}, there exists a vector in $\Sigma_0$ which is a face of $\Gamma$. 
    Since $\Sigma \subseteq \Gamma$, such a vector cannot be surrounded. Next, suppose that $\Gamma = \bbr^2$. By proposition~\ref{prop:cara1}, there must exist a minimal cone $G=\la a_i,a'_i\ra$ such that $G \cap \Sigma^c$ is non-empty and, by proposition~\ref{prop:a1a1pcov}, $G^\circ \subseteq \Sigma^c$ proving that both $a_i$ and $a'_i$ are not surrounded. 
\end{proof}
 
\begin{proposition}\label{prop:local2}
Suppose $n=2$ and let $i,j$ denote two distinct indices. 
Let $\bar{v}$ denote the orthogonal projection of $v \in \bbr^2$ onto $a_i^{\perp}$.  
Then $a_i$ is surrounded if and only if it is either self or lazily surrounded or $a_i \simeq a_j$ and the pair $\{\ba'_i,\ba'_j\}$ form a $1$-dimensional Q-covering or $a_i \simeq a'_j$, and the pair $\{\ba'_i,\ba_j\}$ form a $1$-dimensional Q-covering. 
\end{proposition}
\begin{proof}
    Sufficiency is immediate using Proposition~\ref{prop:self}. For necessity, assume that $a_i$ is not self nor lazily surrounded. 
    Thus $a_i \in \Sigma_1(a'_i)$. 
    If $a_i \simeq a'_i$, then $\Sigma$ reduces to $\la a_i, a_j\ra \cup \la a_i, a'_j \ra$ and $a_i$ is surrounded if and only if it is self surrounded, contradicting the assumption. So $a_i \not\simeq a'_i$. 
    If $a_i \simeq a_j$ and $a_j \simeq a'_j$ then the surrounding is impossible as $\Sigma$ reduces to one cone having $a_i$ as generator.  
    If $a_j \not\simeq a'_j$, then $a_i$ must be surrounded by $\la a_i, a'_j\ra \cup \la a'_i, a_j \ra$ which is effectively equivalent to checking that $\{\ba'_i,\ba'_j\}$ form a $1$-dimensional Q-covering by Proposition~\ref{prop:self}. The same discussion holds when $a_i \simeq a'_j$ and one needs to check that $\{\ba'_i,\ba_j\}$ form a $1$-dimensional Q-covering. 
\end{proof}
Proposition~\ref{prop:local2} is amenable to an algebraic characterization of the Q-covering problem in $n=2$. 
For lazy surrounding, checking if $a_i \in \la a'_i, a_j \ra^\circ$, $j \neq i$, amounts to simply checking that the determinants of the three  matrices $\begin{pmatrix} a'_i & a_j\end{pmatrix}$, $\begin{pmatrix} a'_i & a_i\end{pmatrix}$, and $\begin{pmatrix} a_i & a_j\end{pmatrix}$ have the same sign. 
Likewise, checking if $a_i$ is self surrounded amounts to verifying that the determinants of $\begin{pmatrix} a_i & a_j \end{pmatrix}$ and $\begin{pmatrix} a_i & a'_j \end{pmatrix}$, $j \neq i$, have opposite signs. 
We observe that this condition should not be confused with $a_i\in \langle a_j, a'_j\rangle^\circ$, which is only a special case ($a_i$ needs not be in the interior of $\la a_j,a'_j \ra$ to be self surrounded).  
As stated in Proposition~\ref{prop:self}, self surrounding can be equivalently checked by projecting on the orthogonal space of $a_i$ and appealing to the following simple fact. 
\begin{theorem}\label{thm:1D}
    Let $a,a' \in \bbr$. 
    The pair $\{a,a'\}$ defines a Q-covering of $\bbr$ if and only if $a a' < 0$, providing thereby a partition for $\bbr$.  
\end{theorem}
\begin{proof}
    The cones generated by $a$ and $a'$ cover, or more precisely partition, $\bbr$ if and only if $a,a'$ are both nonzero and have opposite signs. 
\end{proof}
\begin{remark}
Observe that, for $n=1$, $\Gamma$ and $\Sigma$ coincide and that $a$ is surrounded if and only if $a \neq 0$. Interestingly, Corollary~\ref{coro:char} does not hold for $n=1$ since $a$ and $a'$ can be both (lazily) surrounded (for instance when $a' \simeq a$) and $\Sigma = \Gamma \subset \bbr$. 
Intuitively, when $a'_i$ is lazily surrounded, it is somehow `redundant' with $a_i$ (with $a'_i \simeq a_i$ being the simplest--and perhaps strongest--form of redundancy). So when $a'_i$ is redundant for all $i$, $\Sigma$ is unlikely to be covering. 
\end{remark}
Algorithm~\ref{alg:local2D} outputs the set of conditions required for $a_i$ to be surrounded according to Proposition~\ref{prop:local2}. When $a_i$ is identically zero, all conditions fail as desired ($a_i=0$ cannot be surrounded). Applying the algorithm to the four involved vectors outputs an algebraic characterization for $\Sigma$ to be covering for $n=2$; we will denote it in the sequel by \texttt{QCovering}$[\{a_1,a'_1\},\{a_2,a'_2\}]$. 
In particular, one gets the following characterization for Q-matrices in dimension $2$. 

\begin{algorithm}[t]
  \DontPrintSemicolon
  \LinesNumbered
  \SetAlgoVlined
  \KwData{Two pairs $\{a_1,a'_1\}$ $\{a_2,a'_2\}$ of vectors in $\bbr^2$.} 
   $\{u_1,u_2\} \leftarrow a_i$  \;
   $u^{\perp} \leftarrow \{-u_2,u_1\}$ \Comment{Orthogonal vector if $a_i \neq 0$} \;
   $c_1 \leftarrow (u^{\perp} . a_j)(u^{\perp} . a'_j) < 0$ \Comment{Self surrounding} \label{line:c1alg4}\; 
   $c_2 \leftarrow \det(\begin{smallmatrix}
       a'_i & a_i 
   \end{smallmatrix})\det(\begin{smallmatrix}
       a_i & a_j 
   \end{smallmatrix}) > 0  \lor \det(\begin{smallmatrix}
       a'_i & a_i 
   \end{smallmatrix})\det(\begin{smallmatrix}
       a_i & a'_j 
   \end{smallmatrix}) > 0$ \Comment{Lazy surrounding} \; 
   $c_3 \leftarrow \det(\begin{smallmatrix}
       a_i & a_j 
   \end{smallmatrix})=0 \land a_i . a_j >0 \land (u^{\perp} . a'_i)(u^{\perp} . a'_j) < 0$ \Comment{$a_i \simeq a_j$} \label{line:c5alg4} \; 
   $c'_3 \leftarrow \det(\begin{smallmatrix}
       a_i & a'_j 
   \end{smallmatrix})=0 \land a_i . a'_j >0 \land (u^{\perp} . a'_i)(u^{\perp} . a_j) < 0$ \Comment{$a_i \simeq a'_j$} \label{line:c6alg4} \;
   \KwRet{$c_1 \lor c_2 \lor c_3 \lor c'_3$} 
  \caption{Surrounding of $a_i$ ($n=2$).}\label{alg:local2D}
\end{algorithm}

\begin{theorem}
    \label{thm:Q2}
    The matrix $\bigl( \begin{smallmatrix} m_1 & m_2 \\ m_3 & m_4 \end{smallmatrix} \bigr)$ is a Q-matrix if and only if 
    \begin{equation}\label{eq:qqq}  
\begin{split}
    &(m_1<0\land m_2>0\land m_3>0\land m_4<0 \land m_1 m_4-m_2 m_3<0) \\
    \lor\, &(m_1<0\land m_2>0\land m_3<0\land m_4>0 \land  m_1 m_4 - m_2 m_3>0) \\
    \lor\, &(m_1=0\land m_2>0\land m_3<0\land m_4>0) \\
    \lor\, &(m_1>0\land m_3=0\land m_4>0) \\
    \lor\, &(m_1>0\land m_2\geq 0\land m_4>0) \\
    \lor\, &(m_1>0\land m_2<0\land m_3>0\land m_1 m_4-m_2 m_3>0) \\
    \lor\, &(m_1>0\land m_2<0\land m_3<0\land m_1 m_4-m_2 m_3>0) \enspace . 
\end{split}
\end{equation}
\end{theorem}
For the sake of comparison, we give below the relatively much simpler conditions for $M$ to be a P-matrix requiring that all the principal minors of $M$ to be positive: 
\[
m_1 > 0 \land m_4 > 0 \land m_1 m_4 - m_2 m_3 > 0  \enspace . 
\]
One observes that P-matrices are a special case of Q-matrices since $m_1 > 0 \land m_4 > 0 \land m_1 m_4 - m_2 m_3 > 0$ implies (without being equivalent to) the last four conjunctions of~\eqref{eq:qqq}. 

\begin{remark}\label{rq:Q2sign}
    The fact that Theorem~\ref{thm:Q2} involves only sign conditions on the subdeterminants of the matrix $M$ is not a coincidence. 
    In fact, Algorithm~\ref{alg:local2D} can be equivalently stated in terms of sign conditions of the subdeterminants of the matrix $(\begin{smallmatrix}
        a_1 & a_2 & a'_1 & a'_2
    \end{smallmatrix})$. Indeed, on one hand, the scalar product $u^{\perp} . a_j$ in Line~\ref{line:c1alg4} is equal to $\det(\begin{smallmatrix} a_i & a_j \end{smallmatrix})$. The same holds for the other scalar products involving $u^{\perp}$. 
    On the other hand, for $v,w \in \bbr^2$, the equivalence $v \simeq w$ characterized by $\det(\begin{smallmatrix}
       v & w 
   \end{smallmatrix})=0 \land v. w >0$ in Algorithm~\ref{alg:local2D}, can be reformulated as  
    \begin{equation}\label{eq:subsub}
   (\det(\begin{smallmatrix}
       v & w 
   \end{smallmatrix})=0) \land (v_1 w_1 > 0 \lor v_2 w_2 > 0), 
    \end{equation} 
   making explicit the subdeterminants. Therefore, the conditions $c_3$ (Line~\ref{line:c5alg4}) and $c'_3$ (Line~\ref{line:c6alg4}) are also amenable to sign conditions on appropriate subdeterminants.  
\end{remark}

To get an algebraic characterization for self surrounding for $n=3$, in addition to the aforementioned characterization of the Q-covering problem for $n=2$, we further need means to perform the projection on the orthogonal space of a vector $u=(u_1,u_2,u_3) \neq 0$. 
To do so, we use one of the following generic projectors: 
\begin{equation}\label{eq:genericpi}
    \pi_u=\begin{pmatrix}
         u_2 & -u_1 & 0 \\
         u_3 & 0 & -u_1
     \end{pmatrix}, \text{ or }  
     \enspace  
     \pi_u=\begin{pmatrix}
         -u_2 & u_1 & 0 \\
         0 & u_3 & -u_2
     \end{pmatrix}, \text{ or }  
     \enspace 
     \pi_u=\begin{pmatrix}
         0 & -u_3 & u_2 \\ 
         -u_3 & 0 & u_1 
     \end{pmatrix} \enspace .
\end{equation}
%
Setting $u$ to $a_i$, one gets that $a_i$ is self surrounded if and only if $\bbr^2 \subseteq \{\pi_{u} a_j, \pi_{u} a'_j\}\oplus\{\pi_{u} a_k, \pi_{u} a'_k\}$ as shown in Algorithm~\ref{alg:local3D}. Notice that, when $u=0$, $\pi_u=0$ and the algorithm returns False as expected (the vector $0$ cannot be self surrounded). 

\begin{lemma}\label{lem:self3sub}
    Let $u,v,w \in \bbr^3$ and let $\pi_u$ denote a generic projector (cf. Eq.~\eqref{eq:genericpi}). 
    Deciding the sign of $\det(\begin{smallmatrix}
        \pi_u v &\ \pi_u w
  \end{smallmatrix})$ and the equivalence $\pi_u v \simeq \pi_u w$ reduce to sign conditions on the subdeterminants of the matrix $(\begin{smallmatrix}
      u & v & w
  \end{smallmatrix})$. 
\end{lemma}
\begin{proof}
    Suppose (a similar discussion holds for the other projectors) 
    \[
    \pi_u=\begin{pmatrix}
         u_2 & -u_1 & 0 \\
         u_3 & 0 & -u_1
     \end{pmatrix}, \quad  
     \begin{array}{l l}
          s_1 = u_2 v_1 - u_1 v_2  \quad   &  s_3 = u_2 w_1 - u_1 w_2  \\ 
          s_2 = u_3 v_1 - u_1 v_3  \quad   &  s_4 = u_3 w_1 - u_1 w_3 
     \end{array}
    \] 
    Then, one has  
\[ 
\pi_u v = \begin{pmatrix}
s_1 \\ 
s_2     
\end{pmatrix}, \quad \pi_u w = \begin{pmatrix}
    s_3 \\
    s_4
\end{pmatrix}, \quad
  \det(\begin{smallmatrix}
        \pi_u v & \pi_u w
  \end{smallmatrix}) 
  = u_1 \det(\begin{smallmatrix}
        u & v & w 
  \end{smallmatrix}),  
\]
Thus, as mentioned in Remark~\ref{rq:Q2sign}, Eq.~\eqref{eq:subsub}, the condition  
\[
\det(\begin{smallmatrix}
       \pi_u v &\ \pi_u w 
   \end{smallmatrix})=0 \ \land \ (\pi_u v) . (\pi_u w) > 0, 
\]
becomes equivalent to 
\[
u_1 \det(\begin{smallmatrix}
        u & v & w 
  \end{smallmatrix}) = 0  \ \land \ 
(s_1 s_3 > 0 \lor s_2 s_4 > 0), 
\]
making therefore explicit the sign conditions on the subdeterminants of $(\begin{smallmatrix}
    u & v & w
\end{smallmatrix})$.   
\end{proof}  
As already observed, self surrounding for $n=3$ reduces to a planar Q-covering problem (Proposition~\ref{prop:self}) which is in turn equivalent to four surrounding problems for $n=2$ (Theorem~\ref{thm:n2surrounding}), each characterized in Proposition~\ref{prop:local2}, and implemented in Algorithm~\ref{alg:local2D}.  
Lemma~\ref{lem:self3sub} is the last ingredient to show that self surrounding for $n=3$ reduces to sign conditions on the subdeterminants of the matrix $(\begin{smallmatrix}
        a_1 & a_2 & a_3 & a'_1 & a'_2 & a'_3
    \end{smallmatrix})$. 

\begin{algorithm}[t]
  \DontPrintSemicolon
  \LinesNumbered
  \SetAlgoVlined
  \KwData{Three pairs $\{a_1,a'_1\}$ $\{a_2,a'_2\}$ $\{a_3,a'_3\}$ of vectors in $\bbr^3$.} 
   $\{u_1,u_2,u_3\} \leftarrow a_i$  \; 
   $c \leftarrow \pi_u = \bigl(\begin{smallmatrix}
         u_2 & -u_1 & 0 \\
         u_3 & 0 & -u_1
     \end{smallmatrix}\bigr)\ \lor\ \pi_u=\Bigl(\begin{smallmatrix}
         u_2 & -u_1 & 0 \\ 
         0 & u_3 & -u_2 
     \end{smallmatrix}\Bigr)\ \lor\ \pi_u=\Bigl(\begin{smallmatrix}
         0 & u_3 & -u_2 \\ 
         u_3 & 0 & -u_1 
     \end{smallmatrix}\Bigr)$ \; 
   \KwRet{$c \land \texttt{QCovering}[\{\pi_u a_j, \pi_u a'_j\},\{\pi_u a_k, \pi_u a'_k\}]$} 
  \caption{\texttt{SelfSurrounding}: self surrounding of $a_i$ ($n=3$).}\label{alg:local3D}
\end{algorithm}

\subsection{Lazy Covering} 
To get an algebraic characterization for the Q-covering problem for $n=3$, we still need to translate the condition $a_i \in \Sigma(a'_i)$ into an equivalent explicit set of constraints on the involved vectors. 
This task reduces to checking whether a vector belongs to a cone spanned by three vectors. 
First, the equivalence relation $u \simeq v$ is encoded as $u \times v =0 \land u.v >0$ where $u \times v$ denotes the cross product. 
To check whether $u \in \la a_1,a_2,a_3\ra^\circ$, we simply verify that 
$\det\begin{pmatrix}
    a_1 & a_2 & a_3
\end{pmatrix}, \det\begin{pmatrix}
    u & a_2 & a_3
\end{pmatrix}, \det\begin{pmatrix}
    a_1 & u & a_3
\end{pmatrix}, \det\begin{pmatrix}
    a_1 & a_2 & u
\end{pmatrix}$ are all positive or all negative (the ordering of the column vectors of the matrices ensure the coherence of the orientation; cf. the discussion at the beginning of Section~\ref{sec:surrounding}); We can alternatively regard lazy surrounding as a Q-covering problem as stated in Proposition~\ref{lem:lazy}.  
It remains to check whether $u$ belongs to the faces of the cone, namely cones of the form $\la v,w \ra$. 
We do so by checking whether $u$ is equivalent to the generators, $v,w$, or it belongs to the relative interior of $\la v,w\ra$ as detailed in Algorithm~\ref{alg:lazycovering}.   
Finally, Algorithm~\ref{alg:incone} returns an equivalent characterization for $a'_i$ to be in a c-cone rooted at $a_i$. 

\begin{remark}\label{rq:lazy3sub} 
    The conditions of Algorithm~\ref{alg:lazycovering} are also amenable to equivalent conditions on the sign of the subdeterminants of the matrix $(\begin{smallmatrix}
        u & v & w
    \end{smallmatrix})$. 
    For $u \simeq v$, the components of the cross product are subdeterminants by definition. Moreover, the condition $u \times v = 0 \land u . v > 0$ is equivalent to $u \times v = 0 \land (u_1 v_1 > 0 \lor u_2 v_2 > 0 \lor u_3 v_3 > 0)$, making explicit the subdeterminants. For condition $c_3$ (Line~\ref{line:c3alg6}), 
    one first observes that $\det(\begin{smallmatrix}
       v & w & v \times w
   \end{smallmatrix})>0$ is equivalent to $v \times w \neq 0$.~\footnote{Indeed,  $\det(\begin{smallmatrix}
       v & w & v \times w
   \end{smallmatrix}) = \norm{v \times w}^2$ for any vectors $v,w$. Requiring the strict inequality ensures that the relative interior $\la v,w\ra^\circ$ is not empty.}  
   When in addition $\det(\begin{smallmatrix}
       u & v & w \end{smallmatrix})=0$, $u = \alpha v + \beta w$ for some scalars $\alpha,\beta$, and $\det(\begin{smallmatrix}
       v & u & v \times w
   \end{smallmatrix})= \beta \det(\begin{smallmatrix}
       v & w & v \times w
   \end{smallmatrix})$ (the determinant is multilinear and alternating). 
   Thus, assuming $v \times w \neq 0$, $\det(\begin{smallmatrix}
       v & u & v \times w
   \end{smallmatrix}) > 0$ if and only if $\beta>0$. As $u . v^{\perp} = \beta w . v^{\perp}$, one can then reformulate $\beta > 0$ as $\pi_v u \simeq \pi_v w$, which has been shown in Lemma~\ref{lem:self3sub} to reduce to sign conditions on the subdeterminants of the matrix $(\begin{smallmatrix} u & v & w \end{smallmatrix})$. 
   The same discussion holds for $\det(\begin{smallmatrix}
       u & w & v \times w
   \end{smallmatrix})>0$. Summing up, condition $c_3$ of Algorithm~\ref{alg:lazycovering} becomes equivalent to 
   \(
    \det(\begin{smallmatrix}
       u & v & w \end{smallmatrix})=0 \land v \times w \neq 0 \land \pi_v u \simeq \pi_v w \land \pi_w v \simeq \pi_w u 
   \). 
\end{remark}

\begin{algorithm}[t]
  \DontPrintSemicolon
  \LinesNumbered
  \SetAlgoVlined
  \KwData{Three symbolic vectors of dimension $3$.} 
   $c_1 \leftarrow u \times v = 0 \land u . v >0$ \Comment{$u \simeq v$}\; 
   $c_2 \leftarrow u \times w = 0 \land u . w >0$ \Comment{$u \simeq w$}\;
   $c_3 \leftarrow \det(\begin{smallmatrix}
       u & v & w
   \end{smallmatrix})=0 \land \det(\begin{smallmatrix}
       v & w & v \times w
   \end{smallmatrix})>0 \land \det(\begin{smallmatrix}
       v & u & v \times w
   \end{smallmatrix})>0 \land \det(\begin{smallmatrix}
       u & w & v \times w
   \end{smallmatrix})>0$ \Comment{$u \in \la v,w\ra^\circ$} \label{line:c3alg6}\;
   \KwRet{$(u=0) \lor c_1 \lor c_2 \lor c_3$} 
  \caption{\texttt{InFace}: symbolic characterization of $u \in \la v,w\ra$ ($n=3$).}\label{alg:lazycovering}
\end{algorithm}

\begin{algorithm}[t]
  \DontPrintSemicolon
  \LinesNumbered
  \SetAlgoVlined
  \KwData{Four symbolic vectors of dimension $3$.} 
   $\{d_1,d_2,d_3,d_4\} \leftarrow \{\det(\begin{smallmatrix}
    a_i & a_j & a_k
\end{smallmatrix}), \det(\begin{smallmatrix}
    a'_i & a_j & a_k
\end{smallmatrix}), \det(\begin{smallmatrix}
    a_i & a'_i & a_k
\end{smallmatrix}), \det(\begin{smallmatrix}
    a_i & a_j & a'_i
\end{smallmatrix})\}$ \; 
   $c_1 \leftarrow \bigl(\bigwedge_{i=1}^4 d_i>0\bigr) \lor \bigl(\bigwedge_{i=1}^4 d_i<0\bigr)$ \Comment{Topological interior} \; 
   $c_2 \leftarrow \texttt{InFace}[a'_i,\{a_i,a_j\}] \lor \texttt{InFace}[a'_i,\{a_i,a_k\}] \lor \texttt{InFace}[a'_i,\{a_j,a_k\}]$ \Comment{Faces} \; 
   \KwRet{$(a'_i=0) \lor c_1 \lor c_2$} 
  \caption{\texttt{InCone}: symbolic characterization of $a'_i \in \la a_i,a_j,a_k\ra$ ($n=3$).}\label{alg:incone}
\end{algorithm} 

As an immediate corollary of Remark~\ref{rq:Q2sign}, Lemma~\ref{lem:self3sub}, and Remark~\ref{rq:lazy3sub}, we get the following nontrivial result which, until now, was considered an open problem to the best of our knowledge. 
\begin{theorem}\label{thm:sub3n}
    For $n=3$, the Q-covering is characterized by sign conditions on the subdeterminants of the matrix $(\begin{smallmatrix}
        a_1 & a_2 & a_3 & a'_1 & a'_2 & a'_3
    \end{smallmatrix})$. In particular, deciding if a $3$-by-$3$ matrix $M$ is a Q-matrix reduces to sign conditions on the subdeterminants of $M$ effectively constructed by Algorithms 1--7.  
\end{theorem}

It's worth mentioning, however, that~\cite{Garcia1983} showed that for \emph{super-regular} matrices (cf. Section~\ref{sec:special}), when the \emph{conical degree} of the piecewise linear mapping associated to $M$ is not zero, $M$ is Q-matrix. While the conical degree is determined using the signs of the subdeterminants of $M$, it was unclear whether the signs of subdeterminants were enough for the cases where the conical degree is zero, and more generally, for matrices that are not super-regular (for which the concept of conical degree is not well defined). 

All algorithms were implemented to arrive at an algebraic characterization of the Q-covering problem when $n=3$.~\footnote{We used Mathematica. Notebook available here \url{https://gitlab.inria.fr/kghorbal/qmatrices}.} 
When $\Sigma_0 =\{e_1,e_2,e_3,-M_1,-M_2,-M_3\}$ for a matrix $M=(\begin{smallmatrix}
    M_1 & M_2 & M_3
\end{smallmatrix})$, one gets a list of sign conditions on the subdeterminants of $M$ for $M$ to be a Q-matrix. 
Checking if a given instance is a Q-matrix is thus performed in constant time: it suffices to substitute the values and check the conditions. 
For instance, for 
\begin{equation}\label{eq:M}   
M = \begin{pmatrix}
    m_1 & m_2 & m_3 \\
    m_4 & m_5 & m_6 \\ 
    m_7 & m_8 & m_9
\end{pmatrix}, 
\end{equation}
any specialization of $m_1,\dotsc,m_9$ that satisfies 
\[
     m_1 > 0 \land m_2 = 0 \land m_3 = 0 \land m_5 < 0 \land m_6 > 0 \land m_8 > 0 \land m_9 < 0 \land m_5 m_9 - m_6 m_8 < 0   
\]
is a Q-matrix (which is clearly not a P-matrix as $m_5 < 0$). 
Such a characterization turned out to be very useful to automatically find counter examples to sharpen our intuitions and help answering certain conjectures as we shall see next.

\subsection{Generating Special Q-matrices}\label{sec:special} 
\cite{karamardian1972complementarity} 
drew a specific attention to R$_0$ matrices (also known as super-regular matrices) for which LCP$(0,M)$ has a unique solution. 
It subsequently played an important role in LCP theory. \\
\cite{aganagic1979note} proved later that among P$_0$ matrices (i.e. matrices for which all principal minors are nonnegative), the subclasses of Q-matrices and R$_0$-matrices are equivalent. 

We wanted to know whether a Q-matrix that is not super-regular exists for $n=3$. It was previously known~\cite[p. 177]{KELLY1979175} that Q-matrices with flat c-cones are possible for $n=3$. 
We were thus interested in finding a Q-matrix with non-pointed and non-flat c-cones. To do so we fixed $a'_2$ to $-e_1$ and found that the following conditions (among others) on the subdeterminants of $M$ (cf. Eq.~\eqref{eq:M}) have to hold 
\[
m_1 > 0 \land m_3 < 0 \land m_6 < 0 \land m_7 > 0 \land 
 m_9 < 0 \land m_1 m_9 - m_3 m_7 > 0 \land m_4 m_9 - m_6 m_7  < 0 .
\]
The following instance is a particular case.  

\begin{example}\label{ex:examplen3}
Consider $-M_1=(-2,-4,-3)$, $-M_2=-e_1$, and $-M_3=(1,1,1)$. The c-cones $\la e_1, -M_2, e_3\ra$ and $\la e_1, -M_2, -M_3\ra$ are non-pointed and therefore the following $3$-by-$3$ matrix $M$ is not super-regular. 
It is however a Q-matrix.~\footnote{The fact that the provided matrix is a Q-matrix can be checked independently using for instance a quantifier elimination procedure over the reals, like the Cylindrical Algebraic Decomposition or Gale's algorithm~\cite[p. 4]{AganagicCottle78}.} 
\[
M = \begin{pmatrix}
 M_1 & M_2 & M_3   
\end{pmatrix} 
= \begin{pmatrix}
2 & 1 & -1 \\
4 & 0 & -1 \\ 
3 & 0 & -1
\end{pmatrix} \enspace .
\]
We can also easily check that there are no Q-matrices of the form $(\begin{smallmatrix}
    v & e_1 & e_2 
\end{smallmatrix})$ for any vector $v$ (which would have lead to $4$ degenerate c-cones that are non-pointed and non-flat).  
\end{example}

\cite{MURTY197265} gave the symmetric matrix $M_{\text{Murty}}$ (see below) for which all vectors in $\Sigma_0$ are both self and lazily surrounded. 
Using the algebraic characterization presented in this work, we checked that when all vectors are lazily surrounded, then they are necessarily also self surrounded (the converse isn't true, it suffices to consider any P-matrix). Enforcing lazy surrounding of all vectors in $\Sigma_0$ is unrelated to the symmetry of the matrix $M$. 
It turns out that Murty's example is an instance of the following conjunction: 
\begin{multline}\label{eq:mmurty}
m_1 < 0 \land m_3 > 0 \land m_4 > 0 \land 
 m_5 < 0 \land m_6 > 0 \land m_7 > 0 \land m_8 > 0 \land 
 m_9 < 0 \\ \land m_1 m_5 - m_2 m_4  < 0 \land 
 m_1 m_9 - m_3 m_7  < 0 \land m_5 m_9 - m_6 m_8  < 0  \enspace .  
\end{multline}
We give below several non symmetric instances.  
\begin{example}\label{ex:exampleMurty}
The following two matrices are Q-matrices for which all vectors in $\Sigma_0$ are both self and lazily surrounded (the space is covered twice by the c-cones). They both satisfy the conditions of Eq.~\eqref{eq:mmurty}. 
\[
M_{\text{Murty}} = \begin{pmatrix}
-1 & 2 & 2 \\
 2 & -1 & 2 \\
 2 & 2 & -1 \\
\end{pmatrix},
\qquad 
M = \begin{pmatrix}
-5 & 4 & 3 \\
 2 & -1 & 1 \\
 2 & 2 & -1 \\
\end{pmatrix} \enspace .
\]
Moreover, there are precisely $3$ additional conjunctions, each involving strict sign conditions on the subdeterminants of $M$, such that all vectors in $\Sigma_0$ are both self and lazily surrounded. We give below one instance for each such conjunction (the sign conditions could be retrieved from the examples).   
\[
\begin{pmatrix}
 -7 & 5 & 1 \\
 -6 & 4 & 1 \\
 -8 & 8 & 1 \\
\end{pmatrix}, \qquad
\begin{pmatrix}
 3 & -9 & 1 \\
 4 & -10 & 1 \\
 16 & -16 & 1 \\
\end{pmatrix}, \qquad
\begin{pmatrix}
 7 & 5 & -1 \\
 12 & 2 & -1 \\
 8 & 4 & -1 \\
\end{pmatrix} \enspace .
\]
\end{example}


%
\section*{Conclusion}
We believe that using minimal cones to better understand holes is worth pursuing. An important question with this regard is whether the particular case of almost c-cones is all one needs in dimensions $\geq 4$. 
The approach is also appealing as it doesn't require particular assumptions on degeneracy (but does assume feasibility). It thus provides an interesting geometric alternative to degree theory. 
For $n\leq 3$, one gets in addition an algebraic characterization involving only the signs of the subdeterminants of the involved matrix.  
So far, however, no particular \emph{pattern} emerged from these conditions (unlike the elegant characterization for P-matrices requiring only the positivity of principal minors).  
From a computational standpoint, although getting an algebraic  characterization was relatively easy (the presented algorithms are straightforward to implement), rewriting the so obtained condition into conjunctions on the signs of subdeterminants was computationally involved as many subcases are either redundant or empty.  
It would be really interesting to try to push the same reasoning for dimension $n=4$ to get yet an additional hint about the polynomials involved as well as potential patterns on their sign conditions.  
We do believe that algebraic characterizations are really helpful to sharpen our intuitions and avoid pitfalls by automatically generating instances with the help of a computer. 

\paragraph{Acknowledgments} We would like to sincerely thank Prof. Walter Morris for his valuable comments and several related references he pointed out while reading an earlier draft of this work. We are also in debt to Jean-Charles Gilbert for his thorough reading and the several technical improvements he suggested. We finally thank the anonymous reviewers for their substantial feedback and several suggestions to help improve the clarity and readability of the paper.   

\bibliographystyle{apalike}
\bibliography{references}



\end{document}